\documentclass[a4paper]{amsart}
\usepackage{hyperref}
\newcommand{\expect}{\mathbb E}
\newcommand{\lr}[1]{\left(#1\right)}
\newcommand{\abs}[1]{\left\vert #1 \right\vert}
\newcommand{\norm}[2]{\left\Vert #1 \right\Vert _{#2}}
\newcommand{\set}[1]{\left\{#1\right\}}
\newcommand{\scalar}[2]{\langle #1,#2\rangle}
\newcommand{\real}{\mathbb R}
\newcommand{\nat}{{\mathbb N}}
\newcommand{\tr}[1]{\operatorname{tr}\left[#1\right]}
\newcommand{\myspan}{\operatorname{span}}

\newcommand{\ra}{r_\tau}

\newcommand{\ga}{g_\tau}

\newcommand{\opA}{T^\ast T}
\newcommand{\cN}{{\mathcal{N}}}
\newcommand{\aast}{\tau_\ast}
\newcommand{\xmnd}{f_{m,n}}
\newcommand{\bp}{\lr{Q_n T P_m}^{\dag}}
\newcommand{\jott}{\operatorname{j}}
\newcommand{\Dab}{C_{\alpha,\beta}}
\newcommand{\tra}{t_{R,\alpha}}
\newcommand{\IP}{\mathrm{IP}}
\newcommand{\DP}{\mathrm{DP}}
\newenvironment{frcseries}{\fontfamily{frc}\selectfont}{}
\newcommand{\textfrc}[1]{{\frcseries#1}}
\newcommand{\littleo}{\textfrc o} 

\theoremstyle{plain}
\newtheorem{thm}{Theorem}
\newtheorem{lemma}{Lemma}
\newtheorem{proposition}{Proposition}
\newtheorem{corollary}{Corollary}
\newtheorem{de}{Definition}
\theoremstyle{remark}
\newtheorem{remark}{Remark}
\newtheorem{xmpl}{Example}
\theoremstyle{definition}
\newtheorem{ass}{Assumption}

\numberwithin{equation}{section}
\numberwithin{thm}{section}
\numberwithin{lemma}{section}
\numberwithin{remark}{section}
\numberwithin{proposition}{section}
\numberwithin{corollary}{section}

\title[General regularization schemes for signal detection]{General
  regularization schemes for signal detection  in inverse
  problems} 
\author{Cl\'{e}ment Marteau}
\address{Institut Math\'ematiques de Toulouse, INSA de Toulouse,
  Universit\'e de Toulouse}
\email{clement.marteau@math.univ-toulouse.fr}
\author{Peter Math\'e}
\address{Weierstra{\ss} Institute for Applied Analysis and
  Stochastics, Mohrenstra\ss e 39, 10117 Berlin,  Germany}
\email{peter.mathe@wias-berlin.de}
\date{Version: \today}

\begin{document}
\begin{abstract}
 The authors discuss how general regularization schemes, in particular
 linear regularization schemes and projection schemes, can be used to
 design tests for signal detection in statistical inverse problems. It
 is shown that such tests can attain the minimax separation rates when the
 regularization parameter is chosen appropriately. It is also shown
 how to modify these tests in order to obtain (up to a $\log\log$
 factor) a test which adapts to the unknown smoothness in the
 alternative.
Moreover, the authors discuss how the so-called \emph{direct} and
\emph{indirect} tests are related via interpolation properties. 
\end{abstract}
\maketitle

\section{Introduction and motivation}
\label{sec:intro}

Statistical inverse problems have been intensively studied over the
last years. Mainly, estimation of indirectly observed signals was
considered. On the other hand, there are only a few studies concerned
with signal detection, which is a problem of statistical testing. This
is the core of the present paper. Precisely, we consider a statistical
problem in Hilbert space, where we are given two Hilbert spaces~$H$
and $K$ along with a (compact) linear operator~$T\colon H \to
K$. Given the (unknown) element~$f\in H$ we observe
\begin{equation}
  \label{eq:base}
  Y = T f + \sigma \xi,  
\end{equation}
where $\xi$ is a Gaussian white noise, and $\sigma$ is a positive
noise level. A large amount of attention has been payed to the
estimation issue, where one wants to estimate the function $f$ of
interest, and control the associated error. We refer for instance to
\cite{EHN_1996} for a review of existing methods in a deterministic
setting ($\xi$ is a deterministic error satisfying $\|\xi\|\leq
1$). In the statistical framework, the noise $\xi$ is not assumed to
be bounded.  In this case, there is a slight abuse of notation in
using (\ref{eq:base}). We assume in fact that for all $g\in K$, we can
observe
$$ \scalar Y g = \scalar{T f} g + \sigma \scalar{ \xi} g ,$$
where $\scalar{\cdot}{\cdot}$ denotes the scalar product in $K$.
Details will be given in Section~\ref{sec:construction}.  In this
context, we mention \cite{Bissantz_2007} or \cite{Cavalier_book} among
others for a review of existing methodologies and related rates of
convergence for estimation under Gaussian white noise.

In this study, our aim is to test the null hypothesis that the
(underlying true) signal $f$ corresponds to a given signal $f_{0}$
against a non-parametric alternative. More formally, we test
\begin{equation}
  H_{0}:  f=f_{0}, \ \mathrm{against} \ H_{1,\rho}: \ f - f_{0} \in
  \mathcal E,\ \norm{f - f_{0}}{}\geq \rho,
  \label{eq:model_test}
\end{equation}
where $\mathcal E$ is a subset of $H$, and $\rho>0$ a given
radius. The subset $\mathcal{E}$ can be understood as a smoothness
constraint on the remainder $ f - f_{0}$, while the quantity~$\rho$
measures the amount of signal, different from $f_{0}$, available in
the observation. Following the setting, (\ref{eq:model_test}) is known
as a goodness-of-fit or a signal detection (when $f_0=0$) testing
problem. In the direct case, i.e. when $T=Id$, this problem has been
widely investigated. We mention for instance seminal investigations
proposed in \cite{IngsterI,IngsterII,IngsterIII}.  We refer also to
\cite{MR1935648} where a non-asymptotic approach is proposed.

Concerning testing in inverse problems there exists, up to our knowledge, only  few references, as e.g. \cite{ISS_2012} and
\cite{LLM_2012}. In these contributions, a preliminary estimator $\hat
f$ for the underlying signal $f$ is used. This estimator is based on a
spectral cut-off scheme in \cite{LLM_2012}, or on a refined version
using Pinsker's filter in \cite{ISS_2012}. All these approaches are
based on the same truncated singular value decomposition (see Section
\ref{sec:linreg} for more details). Here we shall consider
general \emph{linear estimators} $\hat f = R Y$, using the data $Y$.
Plainly, since $f_{0}$ and hence $T f_{0}$ are given, we can constrain
the analysis to testing whether $f= 0$ (no signal) against the
alternative $H_{1,\rho}: \ f \in \mathcal E,\ \norm{f}{}\geq \rho$,
and we discuss this simplified
model from now on. \\

In the following, we will deal with level-$\alpha$ tests,
i.e. measurable functions of the data with values in $\lbrace 0, 1
\rbrace$. By convention, we reject $H_0$ if the test is equal to $1$
and do not reject this hypothesis, otherwise. We are interested in the
optimal value of $\rho$ (see (\ref{eq:model_test})) for which a
prescribed level for the second kind error can be attained. More
formally, given a fixed value of $\beta \in ]0,1[$ and a
level-$\alpha$ test, we are interested in the radius
$\rho(\Phi_\alpha,\beta,\mathcal{E})$ defined as
$$
\rho(\Phi_\alpha,\beta,\mathcal{E}) = \inf \left\lbrace \rho \in
  \real^+: \ \sup_{f\in \mathcal{E}, \ \|f\|>\rho}
  P_f(\Phi_\alpha =0) \leq \beta \right\rbrace.
$$
From this, the minimax separation radius
$\rho(\alpha,\beta,\mathcal{E})$ can be defined as the smallest radius
over all possible testing procedures, i.e.
$$ 
\rho(\alpha,\beta,\mathcal{E}) = \arg \min_{\Phi_\alpha}
\rho(\Phi_\alpha,\beta,\mathcal{E}),
$$
and the minimum is over all level-$\alpha$ tests $\Phi_{\alpha}$.  We
stress that this minimax separation radius will depend on the noise
level~$\sigma$, and on spectral properties, both of the operator~$T$
which governs the equation~(\ref{eq:base}), and of the class $\mathcal
E$, describing the smoothness of the alternative.

Lower (and upper) bounds have already been established in order to
characterize the behavior of this radius for different kind of
smoothness assumptions (see for instance \cite{ISS_2012} or
\cite{LLM_2012}). Recent analysis of (classical) inverse problems
adopts a different approach by measuring the smoothness inherent in
the class $\mathcal E$ \emph{relative} to the operator~$T$. By doing
so, a unified treatment of moderately, severely and mildly ill-posed
problems is possible. We take this paradigm here and consider the
classes ~$\mathcal E$ as \emph{source sets}, see details in
\S~\ref{sec:sepradius}.

Also, previous analysis was restricted to the truncated singular value
decomposition of the underlying operator~$T$. This limits the
applicability of the test procedures, since often a singular value
decomposition is hardly available, for instance when considering
partial differential equations on domains with noisy boundary data.
Therefore, the objective in this study is to propose alternative
testing procedures that match the
previous minimax bounds.\\

To this end we first consider general linear regularization in terms
of an operator~$R$ (Sections~\ref{sec:1stkind} \& \ref{sec:2ndkind}),
and we shall then specify these as linear regularization (in
Section~\ref{sec:linreg}) or projection schemes (in
Section~\ref{sec:projschemes}), respectively. In each case, we derive
the corresponding minimax separation radii. Next the relation between
testing based on the estimation of $f$ (inverse test), and test based
on the estimation of $T f$ (direct test) is discussed in Section
\ref{sec:d+i}. Such discussion can already be found
in~\cite{MR2763215}. However, here we highlight that the relation
between both problems can be seen as a result of \emph{interpolation}
between smoothness spaces, the one which describes the signal $f$ and
the one which characterizes the smoothness of $T f$.

Finally, we shall establish in Section \ref{s:adaptation} an adaptive
test, which is based on a finite family of non-adaptive tests.  It
will be shown that this adaptive test, with an appropriately
constructed finite family, is (up to a $\log\log$ factor) as good as
the best among the whole family of tests.

\section{Construction and calibration of the test}
\label{sec:construction}
Considering the testing problem~(\ref{eq:model_test}), most of the
related tests are based on an estimation of $\| f \|^2$ ($\| f-f_0
\|^2$ in the general case). Then, the idea is to reject $H_0$ as soon
as this estimation becomes too large with respect to a prescribed
threshold. As outlined above, in order to estimate $\| f \|^2$ where
$f\in H$, from the observations $Y$, cf. (\ref{eq:base}), we shall use
a general linear reconstruction operator~$R\colon K \to H$. 

\subsection{Notation and assumptions}
First we will specify the assumption on the noise~$\xi$
in~(\ref{eq:base}). 
\begin{ass}[Gaussian white noise]\label{ass:GWN}
The noise $\xi$ is a weak random element  in~$K$, which has absolute
weak second  moments. Specifically, for all $g,g_1,g_2 \in K$, we have
$$ \langle \xi, g \rangle \sim \mathcal{N}(0,\|g\|^2), \quad ( \mathrm{and}
\ \expect \left[\langle \xi,g_1 \rangle \langle \xi,g_2 \rangle
\right] = \langle g_1,g_2 \rangle).$$ 
\end{ass} 
Notice that the second property is a consequence of the first, because
bilinear forms in Hilbert space are determined by their values at the diagonal.
 Under such assumption, given any linear reconstruction operator~$R\colon K \to
 H$ the element $R Y$ belongs to $H$ almost surely, provided that $R$
 is a Hilbert--Schmidt operator (Sazonov's Theorem).  
When specifying the reconstruction~$R$ in Sections~\ref{sec:linreg} \&
\ref{sec:projschemes}, we shall always make sure that this is the case.
  Then the application of $R$ to the
data $Y$ may be decomposed as
\begin{equation}
  \label{eq:estimator}
  R Y = R T f + \sigma R\xi = f_R + \sigma R\xi,\quad f \in H,
\end{equation}
where $f_{R}:=R T f$ denotes the noiseless (deterministic part) of
$RY$.  Along with the reconstruction $RY$ the following quantities
will prove important.  First, we can compute the bias variance
decomposition
\begin{equation}
  \label{eq:bias-var}
  \expect\norm{RY}{}^{2} = \norm{RT f}{}^{2} + \sigma^{2}
  \expect\norm{R\xi}{}^{2} = \norm{f_{R}}{}^{2} + S_R^{2},
\end{equation}
where we introduce the \emph{variance} of the estimator as
\begin{equation}
  \label{eq:s2}
  S_{R}^{2}:=  \sigma^{2}\expect\norm{R\xi}{}^{2} =  \sigma^{2}\tr{R^{\ast}R},
\end{equation}
which is finite if $R$ is a Hilbert--Schmidt operator. In addition the following
\emph {weak variance} will play a role.
\begin{equation}
  \label{eq:v2}
  v_{R}^{2}:= \sigma^{2}\sup_{\norm{w}{}\leq 1}\expect\abs{\scalar{R\xi}{w}}^{2} =
  \sigma^{2} \norm{R}{}^{2}.
\end{equation}
Below, if $R$ is clear from the context we sometimes abbreviate
$S=S_{R}$ and $v=v_{R}$.
 
We will need more precise representation of the trace and norm as
above in terms of the representation of the operator~$R$. Suppose that
we have given $R$ in terms of its singular value decomposition as
\begin{equation}
  \label{eq:R-repr}
  R g = \sum_{j=1}^{\infty} \lambda_{j} \scalar{\psi_{j}}{g} \phi_{j},\quad
  g\in K, 
\end{equation}
where we assume that both sequences $\set{\psi_{j}}_{j\in\nat}$ and
$\set{\phi_{j}}_{j\in\nat}$ are orthonormal bases in $K$ and $H$,
respectively. Moreover, the sequence $\lambda_{j},\ j=1,2,\dots$ is
assumed non-negative and arranged in non-increasing order.  Then the
following is well-known.
\begin{lemma}
  Let $R$ be as in~(\ref{eq:R-repr}). Then
  \begin{enumerate}
  \item $\tr{R^{\ast}R} = \sum_{j=1}^{\infty} \lambda_{j}^{2}$, and
  \item $\norm{R}{}^{2} = \sup_{j=1}^{\infty}\lambda_{j}^{2}$.
  \end{enumerate}
\end{lemma}
From this we can see that~ $v_{R}^{2}\leq S_{R}^{2}$, and typically
these quantities
differ by order. Some explicit computations will be provided below.\\

\subsection{Construction of the test and control of the first kind error}
\label{sec:1stkind}

We see from~(\ref{eq:bias-var}) that the quantity
\begin{math}
  \norm{RY}{}^{2}- S_{R}^{2}
\end{math}
is an unbiased estimator for the norm of $\norm{f_{R}}{}^{2}$. If $R$ is chosen appropriately, this
term is an approximation of $\| f\|^2$, whose value is of first
importance when considering the problem
(\ref{eq:model_test}). Therefore, we shall use a threshold for
$\norm{RY}{}^{2}- S_{R}^{2}$ to describe the test.

Let $\alpha \in (0,1)$ be the prescribed level for the first kind
error, and we agree to abbreviate $x_{\alpha}:= \log(1/\alpha)$. We
define the test $\Phi_{\alpha,R}$ as 
\begin{equation}
\Phi_{\alpha,R} = \mathbf{1}_{\left\lbrace \| RY \|^2 - S_{R}^{2} > \tra \right\rbrace},
\label{def_test}
\end{equation}
where $\tra$ denotes the $1-\alpha$ quantile of the variable
$\| RY\|^2 - S_{R}^{2}$ under $H_0$. Due to the definition of the
threshold~$\tra$, the test~$\Phi_{\alpha,R}$ is a level-$\alpha$ test. Indeed
$$ P_{H_0}(\Phi_{\alpha,R}=1) = P_{H_0}(\| RY \|^2 - S_{R}^{2} > \tra)=\alpha.$$
We emphasize that under $H_{0}$ the
distribution of $\| RY\|^2 - S_{R}^{2} = \sigma^{2}(\norm{R\xi}{}^{2}
-\tr{R^{\ast}R})$ only depends on the chosen reconstruction $R$. Hence
the quantile can be determined, at least approximately.  Proposition
\ref{minoration_quantile} below establishes an upper  bound for this
quantile.  

\begin{proposition}
\label{minoration_quantile}
Let $\alpha$ be a fixed level. Then
$$ \tra \leq 2 \sqrt{ 2x_{\alpha}}  S_{R}v_{R}
+ 2 v_{R}^{2}x_{\alpha},
$$ 
where the quantities $S_{R}^{2}$ and $v_{R}^{2}$ have been introduced
in~(\ref{eq:s2}) and~(\ref{eq:v2}).
\end{proposition}
\begin{proof}
 First notice that under $H_{0}$, $\| RY\|^2= \|  \sigma R\xi
  \|^2$. Then we get
\begin{eqnarray*}
\lefteqn{P_{H_0}(\|RY\|^2 -S_R^2 >2 \sqrt{ 2x_{\alpha}}  S_{R}v_{R}
+ 2 v_{R}^{2}x_{\alpha})}\\
&  =  & P_{H_0} \lr{ \| \sigma R\xi \|^2 - S_R^{2} > 2 \sqrt{ 2x_{\alpha}}  S_{R}v_{R}
+ 2 v_{R}^{2}x_{\alpha}}\\
& =   & P_{H_0} \lr{ \|  \sigma R\xi \|^2 - \expect\|  \sigma R\xi
  \|^2 > 2 \sqrt{ 2x_{\alpha}} S_R v_R + 2 v_R^2 x_{\alpha}}\\
& \leq & \exp \lr{ - \frac{2x_{\alpha} v_R^2}{2v_R^2}} = \alpha,
\end{eqnarray*}
where we have used Lemma \ref{TL} with 
$ x = \sqrt{2x_{\alpha}}v_R,$
in order to get the last inequality. Hence,
$$ P_{H_0}(\|RY\|^2 -S_R^2 >2 \sqrt{ 2x_{\alpha}}  S_{R}v_{R}
+ 2 v_{R}^{2}x_{\alpha}) \leq \alpha,$$
which leads to the desired result.
\end{proof}

\subsection{Controlling the second kind error}
\label{sec:2ndkind}
Here, our aim is to control the second kind error
by some prescribed level $\beta>0$, and again we
abbreviate~$x_{\beta}:= \log(1/\beta)$. To this end, we have to exhibit
conditions  on $f$ for which the probability
$P_{f}(\Phi_{\alpha} = 0)$ will be bounded by $\beta$.
By construction of the above test this amounts to bounding
\begin{eqnarray}
 P_{f}(\Phi_{\alpha} = 0) 
& = & P_{f}(\norm{RY}{}^{2} - S^{2}\leq \tra) \nonumber \\
& = & P_{f}(\norm{RY}{}^{2} - \expect\norm{RY}{}^{2} \leq  \tra +
S^{2} -  \expect\norm{RY}{}^{2}) \nonumber  \\
& = & P_{f}(\norm{RY}{}^{2} - \expect\norm{RY}{}^{2} \leq \tra
- \norm{f_{R}}{}^{2}),
\label{eq:2kind_error}
\end{eqnarray}
where the latter follows from~(\ref{eq:bias-var}). In this section, we will investigate the lowest possible value of $\| f_R \|^2$ for which the previous probability can be bounded by $\beta$. 

Let $\beta\in ]0,1[$ be fixed. For all $f\in H$, we denote by $t_{R,\beta}(f)$ the $\beta$-quantile of the variable $\|RY \|^2 -S^2$. In other words
\begin{equation}
P_f \left( \|RY \|^2 -\expect \|RY\|^2 \leq t_{R,\beta}(f) \right) = \beta.
\label{eq:beta_quantile}
\end{equation} 
Then, we get from (\ref{eq:2kind_error}) and (\ref{eq:beta_quantile}) that $P_f(\Phi_{\alpha,R}=0)$ will be bounded by $\beta$ as soon as
\begin{equation}
t_{R,\alpha} - \|f_R\|^2 \leq t_{R,\beta}(f) \Leftrightarrow  \|f_R\|^2 \geq t_{R,\alpha} - t_{R,\beta}(f).
\label{eq:condition}
\end{equation}
We have already an upper bound on the $1-\alpha$-quantile
$t_{\alpha,R}$. In order to conclude this discussion, we need a lower
bound on~$t_{R,\beta}(f)$.
\begin{lemma}\label{lem:minoration_quantile1}
  Let the reconstruction~$R$ be given as in~(\ref{eq:R-repr}), and let
\begin{equation}
  \label{eq:sigma}
  \Sigma := \sum_{j=1}^{\infty} \sigma_{j}^{4} + 2 \sum_{j=1}^{\infty}
  \sigma_{j}^{2} \theta_{j}^{2}.
\end{equation}
Then 
\begin{equation*}
t_{R,\beta}(f)  \geq - 2 \sqrt{\Sigma x_{\beta}}.
\end{equation*}
\end{lemma}
\begin{proof}
We first show
the relation of the problem to a specific  sequence space model. By
construction of $R$, using~(\ref{eq:R-repr}), we can expand 
\begin{eqnarray*}
RY = \sum_{j=1}^{\infty} \lambda_{j} \scalar{\psi_{j}}{Y}\phi_{j}
& = & \sum_{j=1}^{\infty} \lambda_{j} \scalar{\psi_{j}}{T f}\phi_{j}
+ \sigma \sum_{j=1}^{\infty} \lambda_{j} \scalar{\psi_{j}}{\xi}\phi_{j},\\
& = & \sum_{j=1}^{\infty} \theta_j \phi_{j}
+  \sum_{j=1}^{\infty} \sigma_j \varepsilon_j \phi_{j},
\end{eqnarray*}
where $ \theta_{j}:=\lambda_{j}
\scalar{\psi_{j}}{T f}$ and $\sigma_{j}:= \sigma \lambda_{j}$ for all
$j \in \nat$, and the $\varepsilon_j$ are i.i.d. standard Gaussian
random variables. Then, we can apply Lemma~\ref{lem:deviation_left}, which gives
$$
P(\norm{RY}{}^{2} - \expect\norm{RY}{}^{2} \leq - 2
\sqrt{\Sigma x_{\beta}}) \leq \beta,
$$  
which completed the proof.
\end{proof}

We are now able to find a condition on $\|f_R\|^2$ in order to control
the second kind error. We introduce the following quantity
\begin{equation}
  \label{eq:dabast}
  C_{\alpha,\beta}^{\ast}=(4\sqrt{x_\beta} + 4\sqrt{2x_\alpha}),
\end{equation}
which is a function of $\alpha$ and $\beta$, only.
\begin{proposition}
\label{pro:cond_fR}
Let us consider the test $\Phi_{\alpha,R}$ as introduced
in~(\ref{def_test}), and let
\begin{equation}
  \label{eq:r-def}
 r^{2}(\Phi_{\alpha,R},\beta):=C_{\alpha,\beta}^{\ast}Sv +  (4x_\alpha + 8x_\beta)v^2 . 
\end{equation}
Then 
$$ 
\sup_{f, \|f_R \|^2 \geq r^{2}(\Phi_{\alpha,R},\beta)} P_f(\Phi_{\alpha,R}
=0) \leq \beta.
$$
\end{proposition}

\begin{proof}
The equation (\ref{eq:condition}) provides a condition for which $P_f(\Phi_\alpha=0) \leq \beta$. Using Proposition \ref{minoration_quantile} and Lemma~\ref{lem:minoration_quantile1}, we see that this condition is satisfied as soon as
$$ \| f_R\|^2 \geq 2\sqrt{\Sigma x_\beta} + 2 \sqrt{2x_\alpha} Sv +2v^2 x_\alpha.   $$
Now we bound
\begin{eqnarray*}
 \Sigma & = & \sigma^4 \sum_{j=1}^{+\infty} \lambda_j^4 + 2 \sigma^2 \sum_{j=1}^{+\infty} \lambda_j^2 \times \lambda_j^2 \langle \psi_j,T f \rangle^2,\\
& \leq & S^2 v^2 + 2 v^2 \| f_R\|^2.
\end{eqnarray*}
Using the inequality ($ab \leq a^2/2 + b^2/2$ for all $a,b\in \real$), we get
\begin{eqnarray*}
2\sqrt{\Sigma x_\beta} 
& \leq & 2Sv \sqrt{x_\beta} + 2 \sqrt{2x_\beta} \|f \| v ,\\
& \leq & 2Sv \sqrt{x_\beta} + \frac{1}{2} \| f \|^2 + 4 x_\beta v^2.
\end{eqnarray*}
In particular, the condition (\ref{eq:condition}) will be satisfied as soon as
\begin{eqnarray*}
&  & \frac{1}{2}  \| f_R\|^2 \geq (2\sqrt{x_\beta} + 2\sqrt{2x_\alpha}) Sv + v^2 (2x_\alpha + 4x_\beta).
\end{eqnarray*}
\end{proof}

\begin{remark}\label{rem:Dab}
Please note that the condition on $\| f_R \|^2$ is (as most of the results presented below) non-asymptotic, i.e. we do not require that $\sigma \rightarrow 0$. Using, the property $v\leq S$, we can obtain the simple bound
\begin{equation}
r^{2}(\Phi_{\alpha},\beta) \leq  C_{\alpha,\beta} S v, \ \mathrm{where} \ C_{\alpha,\beta} = 4\sqrt{x_\beta} + 4\sqrt{2x_\alpha} +4x_\alpha + 8x_\beta .
\label{eq:Dab}
\end{equation}
In an asymptotic setting, the value of the constant $C_{\alpha,\beta}$
may  sometimes be improved. 
In particular, the majorization $v\leq S$ is rather rough. In many
cases, we will only deal with the constant $\Dab^{\ast}$, and we refer
to Corollary~\ref{cor:cabast}. 
\end{remark}

\section{Determining the separation radius under smoothness}
\label{sec:sepradius}

We have seen in the previous section that we need to have that
$\norm{f_{R}}{}^{2} \geq C_{\alpha,\beta} S v$  in order to control
the second kind error. Nevertheless, the alternative in
(\ref{eq:model_test}) is expressed in term of a lower bound on
$\|f\|^2$. In this section, we take advantage on the smoothness of $f$
in order to propose a upper bound on the separation radius.   


Using a triangle inequality, we obtain
$$ \| f_{R}\| \geq \|f \| - \| f-f_{R} \|.$$
Hence, $\| f_{R}\|^2 \geq r^{2}(\Phi_\alpha,\beta)$ as soon as
\begin{eqnarray*}
  & & \|f \| - \| f-f_{R} \| \geq r(\Phi_\alpha,\beta),\\
  & \Leftrightarrow & \|f \|^2  \geq \left( r(\Phi_\alpha,\beta) +\| f-f_{R} \| \right)^2,\\
  & \Leftarrow & \|f \|^2  \geq 2 r^2(\Phi_\alpha,\beta) +2 \| f-f_{R} \|^2,
\end{eqnarray*}
In other words, we get from Proposition  \ref{pro:cond_fR} that
\begin{equation}
\sup_{f, \| f \|^2 \geq 2 r^2(\Phi_\alpha,\beta) +2 \| f-f_{R} \|^2} P_f (\Phi_{\alpha,R}=0) \leq \beta.
\label{eq:cond_f}
\end{equation}
Hence we need to make the lower bound on $\|f\|$ as small as
possible. We aim at finding sharp upper bounds for
\begin{equation}
  \label{eq:minR}
  \inf_{R \in \mathcal{R}} \lr{r^{2}(\Phi_{\alpha},\beta) + \norm{f - f_{R}}{}^{2}},
\end{equation}
where the reconstructions $R$ belong to certain families $\mathcal{R}$.
We shall establish order optimal bounds in two cases, the case of
linear regularization and by using projection schemes.

As already mentioned, we shall measure the smoothness relative to the
operator~$T$, and this is done as follows. Since the operator~$T$ is
compact so is the self-adjoint companion $\opA$. 
 The range of $\opA$ is a (dense) subset in $H$, and one may consider
 an element $f$ smooth, if it is in the range of $\opA$.  To be more
 flexible, we shall do this for more general (operator)
 functions~$\varphi(\opA)$. The corresponding operator~$\varphi(\opA)$
 is compact, whenever, $\varphi(t)\to 0$ as $t\to 0$. Therefore, we
 shall restrict to functions with this property. Precisely, we let
\begin{equation}
  \mathcal{E}_\varphi = \set{ h \in H, \ h = \varphi(\opA)
    \omega, \ \mbox{for some} \ \|\omega \| \leq 1},
  \label{ellips}
\end{equation}
for a continuous non-decreasing function~$\varphi$ which obeys
$\varphi(0)=0$ (index function), be a \emph{general source set}.
It was established in~\cite{MR2384768} that each element in $H$ has
some smoothness of this kind, and hence the present approach is most
general. Examples, which relate Sobolev type balls to the present
setup are given in Examples~\ref{xmpl:moderate} \& \ref{xmpl:severe}.

\subsection{Linear regularization}
\label{sec:linreg}
We recall the notion of linear regularization, see
e.g.~\cite[Definition 2.2]{MR2318806}. Such approaches are rather
popular for estimation purpose. In this section, we describe how
these can be tuned in order to obtain suitable tests.

\begin{de}[linear regularization]\label{de:regularization}
A family of functions
$$
\ga\colon (0,\norm{\opA}{}]\mapsto \real,\
  0<\tau \leq \norm{\opA}{},
$$
is called
  regularization if they are piece-wise continuous in~$\tau$ and the
  following properties hold:
  \begin{enumerate}
  \item\label{it:convergence}{For each~$0< t \leq \norm{\opA}{}$ we
    have that~$\abs{\ra(t)}\to 0$ as~$\tau\to 0$;}

  \item \label{it:gamma1} {There is a constant~$\gamma_{1}$ such
    that~$\sup_{0\le t\le \|\opA\|} \abs{\ra(t)}\leq \gamma_{1}$ for all~$0<\tau \leq \norm{\opA}{}$;}

  \item\label{it:gammaast}{There is a constant~$\gamma_{\ast}\geq 1$ such
    that~$\sup_{0\le t \leq
      \norm{\opA}{}} \tau \abs{\ga(t)}\leq \gamma_{\ast}$ for all
    ${0< \tau < \infty}$,}
  \end{enumerate}
{where~$\ra(t) := 1 - t \ga(t),\ 0\le t \leq
  \norm{\opA}{},$ denotes the residual function.}
\end{de}
Notice, that in contrast to the usual convention we used the symbol
$\tau$ instead of $\alpha$, as the latter is used as control parameter
for the error of the first kind.

Having chosen a specific regularization scheme $\ga$ we assign as
reconstruction the linear mapping
$R_{\tau} := \ga (\opA) T^{\ast}\colon K \to H$. Notice that now, the element
$f_{R}$ is obtained as $f_{R} = f_{\tau}= \ga(\opA)\opA f$.
\begin{xmpl}[truncated svd, spectral cut-off]\label{xmpl:tsvd}
Let $(s_{j},u_{j},v_{j})_{j\in\nat}$ be the singular value
decomposition of the operator~$T$, i.e.,\ we have that
$$
T f = \sum_{j=1}^{\infty} s_{j} \scalar {f}{u_{j}} v_{j},\quad f\in H,
$$
and the singular numbers $s_{1} \geq s_{2}\dots \geq 0$ are arranged
in decreasing order. With this notation we can use the function
$\ga(t):= 1/t,\ t\geq \tau$ and zero else. This means that we
approximate the inverse mapping of $T$ by the finite
expansion~$R_{\tau} Y:=  \sum_{s_{j}^{2}\geq \tau} \frac 1 {s_{j}}\scalar {Y}{v_{j}} u_{j},\ Y\in K,$
  The condition~$s_{j}^{2}\geq \tau$ translates to an upper bound
  $1\leq j \leq D= D(\tau)$.
The element $f_{\tau}$ is then given as $f_{\tau}=  \sum_{j=1}^{D}\scalar {f}{u_{j}} u_{j}$.
\end{xmpl}
\begin{xmpl} [Tikhonov regularization] \label{xmpl:tikhonov}
Another common linear regularization scheme is given with $\ga(t) = 1/(t + \tau),\ t,\tau>0$. In this
  case we have that $R_{\tau} Y = \lr{\tau I + T^{\ast}T}^{-1}
  T^{\ast}Y$, i.e.,\ this is the minimizer of the penalized least
  squares functional $J_{\tau}(f) := \norm{Y - T f}{}^{2} + \tau
  \norm{f}{}^{2},\ f\in H$.
\end{xmpl}
Having chosen any linear regularization, we would like to bound the quantities $S_{\tau}^{2}=S_{R}^{2},v_{\tau}^{2}=v_{R}^{2}$ from~(\ref{eq:s2}),
(\ref{eq:v2}) (with a slight abuse of notation). To this end, we will
impose the following assumption. 
\begin{ass}\label{ass:T-HS} 
  The operator $T$ is a Hilbert-Schmidt operator, i.e.,\ $$ \tr{T^*T} < + \infty.$$ 
\end{ass}
Under the above assumption, the reconstructions~$R_{\tau}$ are also
Hilbert--Schmidt operators, since these are compositions
involving~$T^{\ast}$.

In the following, we shall use the \emph{effective dimension} which
allows to construct a bound on the variance $S_\tau^2$. 
\begin{de}[effective dimension, see~\cite{Cap06,MR2175849}]\label{def:cn}
{The function $\lambda \mapsto \cN(\lambda)$ defined as
\begin{equation}
\label{ass:opreg}
\cN (\lambda) := \tr{(\opA + \lambda I)^{-1}\opA} 
\end{equation}
is called effective dimension of the operator $\opA$ under white noise.}  
\end{de}
By Assumption~\ref{ass:T-HS} the operator $\opA$ has a finite trace,
and the operator $(\opA+\lambda I)^{-1}$ is bounded, thus the function
$\cN$ is finite. 
The following bound is a consequence of~\cite[Lem.~3.1]{B/M2010}.
\begin{equation}
 \tr{\ga^{2}(\opA)\opA} \leq 2 \gamma_{\ast}^{2} \frac{ \cN(\tau)}{\tau},
\label{eq:bound_dim}
\end{equation}
for some constant $\gamma_*>0$.  This, and using the definition of regularization schemes, results in the following bounds.
\begin{lemma}\label{lem:Svbounds}
 Let $R_{\tau} := \ga (\opA) T^{\ast}\colon K \to H$. Assume that Assumption A2 holds, then we have that
\begin{enumerate}
\item[(i)] $\displaystyle{S_{\tau}^{2} \leq 2 \gamma_{\ast}^{2}\sigma^{2}\frac{\cN(\tau)}{\tau},\ \tau>0}$, and
\item[(ii)] $\displaystyle{v_{\tau}^{2} \leq \gamma_{\ast}^{2}\sigma^{2} \frac{1}{\tau},\ \tau>0}$.
 \end{enumerate}
\end{lemma}
\begin{proof}
The proof is a direct consequence of the definition of $S_\tau^2$,
$v_\tau^2$ and of (\ref{eq:bound_dim}).  
\end{proof}

The previous lemma only provides upper bounds for the terms $S_\tau$
and $v_\tau$. For many linear regularization schemes we can actually
show that $v_{\tau} /S_{\tau}\to 0$ as $\tau \to 0$, and we mention the following result.
\begin{lemma}
\label{lem:constants}
 Suppose that the regularization $\ga$ has the following properties.
 \begin{enumerate}
 \item There are constants $\hat c,\hat\gamma >0$ such that
   $\abs{\ga(\hat c \tau)}\geq \hat\gamma/\tau$ for $\alpha>0$, and
 \item for each $0< t \leq \norm{\opA}{}$ the function $\tau \to
   \abs{\ga(t)}$ is decreasing.
 \end{enumerate}
If the singular numbers of the operator~$T$ decay moderately, such
that\\ $\#\set{j,\ \hat c \tau \leq s_{j}^{2} \leq \hat c/{\underline
   c}\tau}\to \infty$ as $\tau\to 0$,  then $\tr{\tau \ga^{2}(\opA)\opA}\to \infty$ as $\tau\to 0$.
Consequently, in this case  we have that $v_{\tau}/S_{\tau}\to 0$ as $\tau\to 0$.
\end{lemma}
\begin{proof}
 For the first assertion we bound, given an $\alpha>0$, and using the
 singular numbers $s_{j}$ of the operator~$T$, the trace as
 follows. We abbreviate, for $s_{j}\geq \hat c \alpha$ the value
 $\beta_{j}:= s_{j}/\hat c$. Then for any $0< \underline c < 1$ we
 find that
 \begin{multline*}
   \tr{\tau \ga^{2}(\opA)\opA}  = \sum_{j=1}^{\infty} \tau
   \ga^{2}(s_{j}^{2})s_{j}^{2}
 \geq \sum_{\underline c s_{j}^{2} \leq \hat c \tau \leq s_{j}^{2}}\tau  \ga^{2}(s_{j}^{2})s_{j}^{2} \\
\geq \sum_{\underline c \beta_{j}\leq \tau \leq \beta_{j}} \tau g_{\beta_{j}}^{2}(\hat c \beta_{j}) \hat
c \beta_{j}
\geq \frac{\lr{\hat\gamma}^{2}}{\hat c}  \sum_{\underline c s_{j}^{2}
 \leq \hat c \tau \leq s_{j}^{2}} \underline c \to \infty,\quad \text{as\ }\tau \to 0.
 \end{multline*}
Finally, by Lemma~\ref{lem:Svbounds} we find that
$$
\frac{v_{\tau}^{2}}{S_{\tau}^{2}} \leq  \gamma_{\ast}^{2} \frac 1
{\tau \tr{\ga^{2}(\opA)\opA}},
$$
and the second assertion is a consequence of the first one.
\end{proof}
\begin{remark}
 The assumptions which are imposed above on $\ga$ are known to hold
 for many regularization schemes, in particular for spectral cut-off
 and (iterated) Tikhonov regularization. The assumption on the
 singular numbers hold for (at most) polynomial decay.
\end{remark} 

Lemma \ref{lem:constants} implies that in some specified cases the
separation radius defined in (\ref{eq:r-def}) is of size
$C_{\alpha,\beta}^\star S_\tau v_\tau$ as $\tau \to 0$. This is summarized in the
following corollary.  

\begin{corollary}\label{cor:cabast}
Let $\Dab^{\ast}$ and  $r^{2}(\Phi_{\alpha,R},\beta)$ be as in~(\ref{eq:dabast})
and~(\ref{eq:r-def}), respectively.
  Under the assumptions of Lemma~\ref{lem:constants} we have that
$$
\frac{r^{2}(\Phi_{\alpha,R},\beta)}{\Dab^{\ast} S_{\tau}v_{\tau}} \to 1 \quad\text{as}\ \tau\to 0.
$$
\end{corollary}

We turn to bounding the bias~$\norm{f -  f_{\tau}}{}$. 
This can be done under the assumption that the
chosen regularization has enough qualification, see e.g.\ \cite{MR2318806}.
\begin{de}[qualification]\label{de:quali}
{Suppose that $\varphi$ is an index function.
The regularization $\ga$ is said to have qualification~$\varphi$ if there is a
constant $\gamma<\infty$ such that
$$
\sup_{0\leq t \leq \| T^*T \|}\abs{\ra(t)}\varphi(t) \leq \gamma \varphi(\tau), \quad \tau>0.
$$}
\end{de}
\begin{remark}
  It is well known that Tikhonov regularization has qualification
  $\varphi(t) = t$ with constant $\gamma=1$, and this is the maximal
  power. On the other hand,   truncated svd has arbitrary
  qualification with constant $\gamma=1$. 
\end{remark}
In this case we can bound the bias at $f_{R}=f_{\tau}$.
\begin{proposition}\label{pro:quali}
  Let $\ga$ be any regularization having qualification~$\varphi$ with
  constant~$\gamma$. If $f\in\mathcal E_{\varphi}$ then
$$
\norm{f - f_{\tau}}{} \leq \gamma \varphi(\tau).
$$ 
\end{proposition}
\begin{proof}
Let $\omega$ with $\|\omega\|\leq1$ such that $f =
\varphi(\opA)\omega$. Then
\begin{displaymath}
\| f_\tau - f \|
 =  \| g_\tau(T^*T) T^*T f -f \| =  \| r_\tau(T^*T) f \| = \| r_\tau(T^*T) \varphi(T^*T) \omega \| \leq \gamma \varphi(\tau).
\end{displaymath}
\end{proof}

Now we have established bounds for all quantities occurring
in~(\ref{eq:minR}), and this yields the main result for linear
regularization.
\begin{thm}\label{thm:main-linear}
Assume that Assumption~\ref{ass:T-HS} holds, and suppose that $\ga$ is a regularization which has
  qualification~$\varphi$, and that~$f\in\mathcal E_{\varphi}$.  Let
  $\aast$ be chosen from the equation
  \begin{equation}
    \label{eq:tauast}
 \varphi^{2}(\tau) = \sigma ^{2} \frac{\sqrt{\cN(\tau)}}{\tau}.    
  \end{equation}
Then, for all $f\in \mathcal{E}_\varphi$, 
$$
\inf_{\tau>0} \lr{r^{2}(\Phi_{\alpha},\beta) + \norm{f - f_{\tau}}{}^{2}} \leq
\lr{C^*_{\alpha,\beta} \sqrt 2 \gamma_{\ast}^{2}+ \frac{(4x_\alpha + 8x_\beta)\gamma_*^2}{\sqrt{\cN(\tau_*)}} + \gamma^{2}}\varphi^{2}(\tau_{\ast}),
$$
where the constant $\Dab^{\ast}$ has been introduced in~(\ref{eq:dabast}). In particular, we get that
$$ \rho^2(\Phi_{\alpha,\tau_*},\beta,\mathcal{E}_\varphi) \leq 2\lr{C^*_{\alpha,\beta} \sqrt 2 \gamma_{\ast}^{2}+ \frac{(4x_\alpha + 8x_\beta)\gamma_*^2}{\sqrt{\cN(\tau_*)}} + \gamma^{2}}\varphi^{2}(\tau_{\ast}).$$
\end{thm}
\begin{proof}
By Proposition~\ref{lem:Svbounds} and Proposition \ref{pro:cond_fR}, we have that
\begin{align*}
\lefteqn{ r^{2}(\Phi_{\alpha},\beta) + \norm{f - f_{\tau}}{}^{2}}\\
 &
=  C^*_{\alpha,\beta} S v + (4x_\alpha + 8x_\beta) v^2 + \norm{f - f_{\tau}}{}^{2},\\
& \leq C^*_{\alpha,\beta} \sqrt 2
\gamma_{\ast}^{2}\sigma^{2}\frac{\sqrt{\cN(\tau)}}{\tau}
+(4x_\alpha + 8x_\beta)\gamma_*^2\sigma^2 \frac{1}{\tau} +  \gamma^{2} \varphi^{2}(\tau),\\
&\leq  \lr{ C^*_{\alpha,\beta} \sqrt 2 \gamma_{\ast}^{2} + \frac{(4x_\alpha + 8x_\beta)\gamma_*^2}{\sqrt{\cN(\tau_*)}} +
  \gamma^{2}}\varphi^{2}(\tau_*), 
\end{align*}
since the parameter $\tau_{\ast}$ equates both terms $\varphi^2(\tau)$ and $\sigma^2 \tau^{-1}\sqrt{\cN(\tau)}$. This 
gives the upper bound.
\end{proof}

\begin{remark}
Up to now, all the presented results are non-asymptotic in the sense that we do not require that $\sigma^2 \rightarrow 0$. In an asymptotic setting, we can remark that $\tau_*$ as defined in (\ref{eq:tauast}) satisfies $\tau_*\rightarrow $ as $\sigma \rightarrow 0$. Since the effective dimension tends to infinity as $\tau \rightarrow 0$, we get that  
$$ \rho^2(\Phi_{\alpha,\tau_*},\beta,\mathcal{E}_\varphi) \leq
2\lr{C^*_{\alpha,\beta} \sqrt 2 \gamma_{\ast}^{2}(1+ \littleo(1))+
  \gamma^{2}}\varphi^{2}(\tau_{\ast}),$$ as $\sigma\rightarrow 0$.
\end{remark}

We shall highlight the above results with two examples. We shall dwell
into these in order to show that the above results are consistent with
other results for inverse testing (see for instance \cite{MR2763215}).

\begin{xmpl}[moderately ill-posed problem]\label{xmpl:moderate}
Let us assume that the singular numbers of the operator~$T$ decay as
$s_{k} \asymp k^{-t},\ k\in\nat$, with $t>1/2$ (in order to ensure
that Assumption~\ref{ass:T-HS} is satisfied).  In this case the effective
dimension asymptotically behaves like~$\cN(\tau) \asymp
\tau^{-1/(2t)}$, as $\tau\to 0$, see for
instance~\cite[Ex.~3]{B/M2010}. The Sobolev ball 
\begin{equation}
  \label{eq:soboloev-ball}
  \mathcal E^{\mathcal X}_{a,2} := \set{f,\ \sum_{j=1}^{\infty}
  a_{j}^{2}\scalar{f}{\phi_{j}^{2}}\leq R^{2}}, \ \mathrm{with} \ a_{j}=j^s,\quad \forall j>1,
\end{equation}
as considered in~\cite{LLM_2012} coincides (up to constants) with
$\mathcal E_{\varphi}$ for the function~$\varphi(u) = u^{s/(2t)},\
u>0$. In this case the value $\tau_{\ast}$ from~(\ref{eq:tauast}) is
computed as $\tau_{\ast} \asymp \sigma^{8t/(4s + 4t + 1)}$, which
results in an asymptotic separation rate of 
$$
 \rho(\Phi_{\alpha,\tau_*},\beta,\mathcal{E}_\varphi) \asymp
 \varphi(\tau_{\ast}) \asymp \sigma^{2s/(2s + 2t + 1/2)},\quad 
\sigma\to 0,
$$
which corresponds to the 'mildly ill-posed case'
in \cite{LLM_2012} or \cite{ISS_2012}, and it is known to be minimax.

\end{xmpl}
\begin{xmpl}[severely ill-posed problem]\label{xmpl:severe}
Here we assume a decay of the form $s_{k} \asymp \exp(-\gamma k),\ k\in\nat$ of the
singular numbers. The effective dimension behaves like $\cN(\tau)\asymp
\frac 1 \gamma \log(1/\tau)$. The Sobolev ball
from~(\ref{eq:soboloev-ball}) is now given as $\mathcal E_{\varphi}$
for a function~$\varphi(u) = \lr{\frac 1 {2\gamma} \log(1/u)}^{-s}
$.
Then the value $\tau_{\ast}$ calculates as $\tau_{\ast}\asymp \sigma^{2}
\lr{\log(1/\sigma^{2})}^{2s +1/2}$, which results in a separation rate
$$
 \rho(\Phi_{\alpha,\tau_*},\beta,\mathcal{E}_\varphi) \asymp \varphi(\tau_{\ast}) \asymp \log^{-s}(1/\sigma^{2}),\quad
\sigma\to 0,
$$
again recovering the corresponding result from~\cite{LLM_2012}.
\end{xmpl}


\subsection{Projection schemes}
\label{sec:projschemes}

Here we follow the ideas from the previous section. Details on the
solution of ill-posed equations by using projection schemes can be
found in~\cite{MR2394505,MR0488721,MR828379}, and our outline follows
the recent~\cite{MR2394505}. In particular we use the
intrinsic requirements such as quasi-optimality and robustness of
projection schemes in order to obtain a  control similar to the
previous section.

 We fix a finite
dimensional subspace~$H_{m}\subset H$, called the \emph{design space}
and/or a finite dimensional 
subspace~$K_{n}\subset K$, called the \emph{data space}. Throughout we
shall denote the 
corresponding orthogonal projections onto~$H_{m}$ by $P_{m}$, and/or
the orthogonal projection onto~$K_{n}$ by~$Q_{n}$. The subscripts $m$
and~$n$ denote the \emph{dimensions} of the spaces. Given such couple~$\lr{H_{m},K_{n}}$ of spaces
we turn from the equation~(\ref{eq:base}) to its discretization
\begin{equation}
  \label{eq:base-discr-eq}
  Q_{n}Y = Q_{n} T P_{m}x + \sigma Q_{n}\xi .
\end{equation}
Without further assumptions, the finite dimensional
equation~(\ref{eq:base-discr-eq}) may have no or many solutions, and
hence we shall turn to the \emph{least-squares solution} as given by
the \emph{Moore-Penrose} inverse, i.e.,  we assign
\begin{equation}
  \label{eq:disc-solu}
  \xmnd := \bp Q_{n}Y.
\end{equation}
\begin{de}[projection scheme, see~\cite{MR2394505}]\label{de:pr-scheme}
{If we are given} 
\begin{enumerate}
\item an increasing sequence~$H_{1} \subset H_{2}\dots \subset H$,
  and 
\item {an increasing sequence~$K_{1} \subset K_{2}\dots \subset K$,
  together with} 
\item {a mapping~$m\to n(m),\ m=1,2,\dots$,} 
\end{enumerate}
{then the corresponding sequence of mappings
\begin{equation}
  \label{eq:bmn-seq}
  Y \to f_{m,n(m)} := \bp Y
\end{equation}
is called~\emph{projection scheme}. }
\end{de}

\begin{xmpl}[truncated svd, spectral cut-off]\label{xmpl:tsvd-proj} The truncated svd, as introduced in
  Example~\ref{xmpl:tsvd} is also an example for a projection scheme,
  if we use the increasing sequences $H_{m} :=
  \myspan\set{u_{1},\dots,u_{m}}\subset H$, and $K_{m}:=
  \myspan\set{v_{1},\dots,v_{m}}\subset K$, respectively. In this case
  we see that $\bp Y = \sum_{j=1}^{m} \frac 1 {s_{j}} \scalar{Y}{v_{j}}$.
\end{xmpl}

Henceforth  we shall
always assume that the mapping $\bp\colon K_{n}\to H_{m}$ is
invertible, i.e.,\ the related linear system of equations has a unique
solution.
This gives an (implicit) relation $n=n(m)$, typically $n=n$ will
do. However, our subsequent analysis will be done using the dimension $m$
of the space $H_{m}$ for quantification. In accordance with this we
will denote $f_{R}$ by $f_{m}$, highlighting the dependence on the dimension.
Thus the linear reconstruction~$R$ is given as $R:= \bp$, and we need
to control $\mathrm{tr}[R^{\ast}R]$ as well as 
$\norm{R}{}$. The latter is related to the robustness (stability) of
the scheme. 
\begin{de}[Robustness]
 { A projection scheme~$\lr{\bp,\ m\in\nat}$ is said to
  be~\emph{robust} if there is a constant~$D_{R}<\infty$ for which
  \begin{equation}
    \label{eq:de-rob}
    \norm{\bp}{} \leq \frac{D_{R}}{\jott(T,H_{m})},\quad m=1,2,\dots
  \end{equation}
Here, the quantity $\jott(T,H_{m})$ denotes the \emph{modulus of
  injectivity} of $T$ with respect to the subspace $H_{m}$, given as 
  \begin{equation}
    \label{eq:jaxm}
    \jott(T,H_{m}) := \inf_{0 \not= x\in H_{m}}\frac{\norm{Tx}{}}{\norm{x}{}}.
  \end{equation}}
\end{de} 
The modulus of injectivity is always smaller than the $m$-th singular
number~$s_{m}= s_{m}(T)$ of the mapping $T$, and hence we say that the subspaces $H_{m}$
satisfy a \emph{Bernstein-type} inequality if there is a constant
$0< C_{B}\leq 1$ such that
$$
C_{B}s_{m}(T) \leq  \jott(T,H_{m}).
$$
We summarize our previous outline as follows.
\begin{lemma}\label{lem:v2bound}
  Suppose that the projection scheme~$\lr{\bp,\ m\in\nat}$ is robust
  and that the spaces $H_{m}$ obey a Bernstein-type inequality. Then
$$
\norm{\bp}{} \leq \frac{D_{R}}{C_{B}} \frac 1 {s_{m}}.
$$ 
In particular we have that 
$$
v_R^2:=v_m^{2} \leq
\sigma^{2}\frac{D_{R}^{2}}{C_{B}^{2}}\frac 1 {s_{m}^{2}}.
$$
\end{lemma}
We turn to bounding $S_R^{2}$. Before doing so we mention that for
spectral cut-off from Example~\ref{xmpl:tsvd-proj},  this bound can easily be established.
\begin{lemma}
For spectral cut-off we have 
$$
S_{R}^{2} = \sigma^{2}\tr{ \lr{\bp}^{\ast} \bp} = 
\sigma^{2}\sum_{j=1}^{m} \frac 1 {s_{j}^{2}}.
$$ 
\end{lemma}
In order to achieve a similar bound in more general situations we need to
impose restrictions on the decay of the singular numbers $s_{j},\
j=1,2,\dots$ 
The use of projection schemes for severely ill-posed problems
requires particular care, and the following restriction, which will be imposed
on the decay of the singular numbers of the operator~$T$ takes this
into account.
We shall assume that the decreasing sequence $s_{j}, j=1,2,\dots,$ is
\emph{regularly varying} for some index $-r$, for some $r\geq 0$, and
we refer to~\cite{MR0333082} for a treatment. In particular this
covers moderately ill-posed problems where $s_{j} \asymp j^{-r}$. We
will not use the index $r$. However, if the sequence $s_{j},
j=1,2,\dots,$ is regularly varying with index $-r$ then the sequence
$s_{j}^{-2}, j=1,2,\dots,$ is regularly varying with index $2r$, and
we have that
$$
\frac 1 m s_{m}^{ 2} \sum_{j=1}^{m} \frac 1 {s_{j}^{2}}
\longrightarrow \frac 1 {2r+1},\quad \text{as}\ m\to\infty. 
$$
In particular there is a constant $C_{r}$ such that 
\begin{equation}
  \label{eq:cr}
  \frac m {s_{m}^{ 2}} \leq C_{r}^{2}  \sum_{j=1}^{m} \frac 1 {s_{j}^{2}},
\end{equation}
and the latter bound is actually all that is needed.
\begin{lemma}\label{lem:S2bound}
  Suppose that the sequence $s_{j}, j=1,2,\dots,$ is such that for the
constant $C_{r}$ the estimate~(\ref{eq:cr}) holds.
If the projection scheme is robust with constant $D_{R}$, and if the
spaces $H_{n}$ obey a Bernstein-type inequality with constant $C_{B}$ then
$$
S_R^2:=S_m^{2} = \sigma^{2} \tr{ \lr{\bp}^{\ast} \bp}
\leq 2 C_{r}^{2} \frac{D_{R}^{2}}{C_{B}^{2}} \sigma^{2}\sum_{j=1}^{m} \frac 1 {s_{j}^{2}}.
$$
If, in addition the Assumption~\ref{ass:T-HS} is satisfied, then we have that
$$
S_R^2:=S_m^{2} \leq C_{r}^{2} \frac{D_{R}^{2}}{C_{B}^{2}}
\sigma^{2}\frac{\cN(s_{m}^{2})}{s_{m}^{2}}. 
$$
\end{lemma}
\begin{proof}
  We notice that the mapping $\lr{\bp}^{\ast}$ is zero
  on   $H_{m}^{\perp}$, the orthogonal complement of $H_{m}$.
So, we take an orthonormal system
$u_{1},u_{2},\dots,u_{m},\dots$, where the first $m$ components span
$H_{m}$.  With respect to this system we see
that
\begin{align*}
 \tr{ \lr{\bp}^{\ast} \bp}& = \tr{ \bp \lr{\bp}^{\ast}}\\
& = \sum_{j=1}^{\infty} \norm{\lr{\bp}^{\ast}  u_{j}}{}^{2}\\
&= \sum_{j=1}^{m} \norm{\lr{\bp}^{\ast}  u_{j}}{}^{2}\\
& \leq m  \norm{\lr{\bp}^{\ast} }{}^{2} = m \norm{\bp }{}^{2}.
\end{align*}
Using Lemma~\ref{lem:v2bound} we see that $ \tr{\lr{\bp}^{\ast}  \bp}\leq 
m \frac{D_{R}^{2}}{C_{B}^{2}}\frac 1 {s_{m}^{2}}.
$
Now we use~(\ref{eq:cr}) to complete the proof of the first assertion.
Under Assumption~\ref{ass:T-HS} we continue and use the inequality $u/v \leq 2 v/(u+v),\ 0< u
\leq v$, to see that 
$$
\sum_{j=1}^{m} \frac 1 {s_{j}^{2}} \leq \frac 2 {s_{m}^{2}}\sum_{j=1}^{m}
\frac{s_{j}^{2}}{s_{j}^{2} + s_{m}^{2}} \leq 2
\frac{\cN(s_{m}^{2})}{s_{m}^{2}}, 
$$
and the proof is complete.  
\end{proof}
\begin{remark}
  Notice that Lemma~\ref{lem:S2bound} provides us with (an order
  optimal) bound for the variance, even if the operator $T$ is not a
  Hilbert--Schmidt one. But, if it is then the obtained bound
  corresponds to the one from Lemma~\ref{lem:Svbounds} (with $\tau
  \leftarrow s_{m}^{2}$).
\end{remark}
Next, we need to bound $\norm{f - f_{R}}{}$, as this was done in
\S~\ref{sec:linreg} by
assuming qualification, and we need a further property of the projection
scheme, called \emph{quasi-optimality}.
We start with the following well-known result, originally from spline
interpolation~\cite{MR0361531}, and used for projection schemes
in~\cite{MR0488721}, which states that
\begin{equation}
  \label{eq:natterer}
  \norm{f -\bp T f }{}\leq \norm{\bp T}{}\norm{f - P_{m}f}{}.
\end{equation}
Therefore, we can bound the bias whenever the norms $\norm{\bp T}{}$
are uniformly bounded.
\begin{de}[quasi-optimality]
 { A projection scheme $Y\to \bp Y$ is \emph{quasi-optimal} if there is
  a constant~$D_{Q}$ such that $\norm{\bp T}{}\leq D_{Q}$.}
\end{de}
We emphasize that under quasi-optimality the bound for the bias
entirely depends on the quality of the projections~$P_{m}$ with respect
to the element~$f$.
\begin{de}[Degree of approximation]
{Suppose that~$\set{H_{m}},\ \dim(H_{m})\leq m$, is a nested set of
design spaces. 
The spaces~$H_{m}$ are said to have the \emph{degree of
  approximation}~$\varphi$
it there is a constant~$C_{D}<\infty$ with
\begin{equation}
  \label{eq:appr-power-r}
\norm{(I  - P_{m})\varphi(\opA)}{} \leq C_{D} \varphi(s_{m+1}),\quad m=1,2,\dots  
\end{equation}}
\end{de}
For
spectral cut-off this bound (with constant $C=1$) is best
possible. Also, using interpolation type inequalities one can verify
this property for many known approximation spaces $H_{m},\
m=1,2,\dots$, we refer to~\cite{MR2394505} for more details on degree
of approximation and Bernstein-type bounds.
We now can state the analogue of Proposition~\ref{pro:quali} for
projection schemes.
\begin{proposition}\label{pro:degree}
  Suppose that the projection scheme is quasi-optimal with
  constant~$D_{Q}$,  and that it has the degree of
  approximation~$\varphi$ with constant~$C_{D}$. If $f\in\mathcal
  E_{\varphi}$ then we have that
$$
\norm{f - f_{m}}{} \leq D_{Q}C_{D} \varphi(s_{m+1}^{2})
$$
\end{proposition}
We now return to the problem raised in~(\ref{eq:minR}). Here, the
family of reconstructions $R$ runs over all projection schemes, and we
can control the bound by a proper choice of the discretization level~$m$.

For the sake of convenience, we will assume in the following that
Assumption~\ref{ass:T-HS} is satisfied, i.e. that $T$ is a Hilbert-Schmidt
operator. If it is not the case, Theorem~\ref{thm:main-projection}
below remains valid when replacing $\sqrt{ \mathcal{N}(s_m^2)}/s_m^2 $
by $\sqrt{\sum_{j=1}^m s_j^{-2}}/s_m$.

\begin{thm}
  \label{thm:main-projection}
Suppose that the approximate solutions are obtained by a projection
scheme which is quasi-optimal and robust and that
Assumption~\ref{ass:T-HS} holds. Furthermore assume that the
design spaces~$H_{m}$ have degree of approximation~$\varphi$ and obey
a Bernstein-type inequality. 
Let $m_{\ast}$ be chosen from
\begin{equation}
  \label{eq:mast-choice}
  m_{\ast}= \max\set{m,\quad \varphi^{2}(s_{m}^{2})   \geq
    \sigma^{2}\frac{\sqrt{\cN(s_{m}^{2})}}{s_{m}^{2} } }. 
\end{equation}

If $f\in\mathcal E_{\varphi}$ then we have that
\begin{eqnarray*}
\lefteqn{\inf_{m}
\set{r^{2}(\Phi_{\alpha},\beta) + \norm{f - f_{m}}{}^{2} }}\\
& \leq &  
\lr{C^*_{\alpha,\beta}\frac{D_{R}^{2}}{C_{B}^{2}}C_{r} + (4x_\alpha + 8x_\beta)  \frac{D_{R}^{2}}{C_{B}^{2}} \frac{1}{\sqrt{\cN(s_{m_*}^{2})}} + 
  D_{Q}^{2}C_{D}^{2}}\varphi^{2}(s_{m_{\ast}}^{2}),
\end{eqnarray*}
where the constant $\Dab^{\ast}$ has been introduced in~(\ref{eq:dabast}).
\end{thm}
\begin{proof}
By using Lemma~\ref{lem:S2bound} and Proposition~\ref{pro:degree} we
see that for any choice of discretization level~$m$ we
have
\begin{align*}
 \lefteqn{r^{2}(\Phi_{\alpha},\beta) + \norm{f - f_{m}}{}^{2}}\\ 
& \leq  C^*_{\alpha,\beta}
\frac{D_{R}^{2}}{C_{B}^{2}}C_{r}\sigma^{2}\frac{\sqrt{\cN(s_{m}^{2})}}{s_{m}^{2}} 
+   (4x_\alpha + 8x_\beta)
\sigma^{2}\frac{D_{R}^{2}}{C_{B}^{2}}\frac 1 {s_{m}^{2}}  +
D_{Q}^{2}C_{D}^{2}\varphi^{2}(s_{m+1}^{2})\\
\begin{split}
\leq  \lr{C^*_{\alpha,\beta}
\frac{D_{R}^{2}}{C_{B}^{2}}C_{r} + (4x_\alpha + 8x_\beta)
\frac{D_{R}^{2}}{C_{B}^{2}} \frac{1}{\sqrt{\cN(s_{m_*}^{2})}}+
D_{Q}^{2}C_{D}^{2}}\times\\
\times\max\set{\sigma^{2}\frac{\sqrt{\cN(s_{m}^{2})}}{s_{m}^{2}},\varphi^{2}(s_{m+1}^{2})}.  
\end{split}
\end{align*}
At the discretization level~$m_{\ast}+1$ we see by monotonicity that
$$
\varphi^{2}(s_{m_{\ast}+1}^{2}) \leq \varphi^{2}(s_{m_{\ast}}^{2}).
$$
Also, by the choice of $m_{\ast}$ we see that
$$
\sigma^{2}\frac{\sqrt{\cN(s_{m_{\ast}}^{2})}}{s_{m_{\ast}}^{2}}\leq
\varphi^{2}(s_{m_{\ast}}^{2}),
$$
hence both terms in the max are dominated by $\varphi^{2}(s_{m_{\ast}}^{2})$,
which allows us to complete the proof.
\end{proof}

Once again, the previous result is non-asymptotic. In the asymptotic
regime, we get the following improvement.
\begin{corollary}
  Under the assumptions of Theorem~\ref{thm:main-projection} we get
  that
\begin{equation*}
\inf_{m}
\set{r^{2}(\Phi_{\alpha},\beta) + \norm{f - f_{m}}{}^{2} }
 \leq  
\lr{C^*_{\alpha,\beta}\frac{D_{R}^{2}}{C_{B}^{2}}C_{r}(1+\littleo(1)) + 
  D_{Q}^{2}C_{D}^{2}}\varphi^{2}(s_{m_{\ast}}^{2}), 
\end{equation*}
as $\sigma\rightarrow 0$.
\end{corollary}
This is an easy consequence of the fact that along with $\sigma\to 0$
we have~$s_{m_{\ast}}^{2}\to 0$, and hence the
effective dimension at $s_{m_{\ast}}^{2}$ tends to infinity. 

\subsection{Discussion}
\label{sec:discussion}

We first highlight the important fact that in both cases (provided
Assumption~\ref{ass:T-HS}  is satisfied), linear
regularization and for projection schemes the upper bound is obtained
by solving the same 'equation', $\sigma^{2} = \tau
\varphi^{2}(\tau)/\sqrt{\cN(\tau)}$, such that relating $\tau_{\ast}
\sim s_{m_{\ast}}^{2}$, see Theorems~\ref{thm:main-linear}
and~\ref{thm:main-projection}. However, this function is different
from the one used for \emph{function estimation} in inverse
problems. In the same setting the 'optimal' parameter $\tau_{est}$ is
there obtained from solving 
$$
 \varphi^{2}(\tau) = \sigma^{2} \frac{\cN(\tau)}{\tau}.
$$ 
Thus, the effective dimension $\cN$, which is designed for estimation
enters in the inverse testing problem in square root, such that
loosely speaking \emph{testing is easier}.

Another remark may be of interest. For the estimation problem, within
the same context, the bias variance decomposition leads to a variance
term $S_{R}^{2}$, and in order to achieve optimal order
reconstruction, this will be calibrated with the function
$\varphi^{2}$. As we have seen above, for testing the same calibration
is done between the functions $S v$ and $\varphi^{2}$. Since, as
already mentioned $S v \leq S^{2}$ this calibration always yields a
smaller value, which again explains the different rates for separation
radius and estimation error.

Previous analysis of the spectral cut-off regularization scheme
for testing in inverse problems revealed the importance of the quantity
\begin{equation}
  \label{eq:rhod}
  \rho_{D} := \lr{ \sum_{j=1}^{D} \frac 1 {s_{j}^{4}}}^{1/4},\quad D=1,2,\dots
\end{equation}
We mention the non-asymptotic lower and upper bounds, slightly adapted to the present
setup, given for instance
in~\cite{LLM_2012} as
\begin{align*}
  \rho^{2}(\mathcal E_{\varphi},\alpha,\beta) & \geq \sup_{D}
  \min\set{c_{\alpha,\beta}^{2}\rho_{D}^{2},\varphi^{2}(s_{D}^{2})},\\
\rho^{2}(\mathcal E_{\varphi},\alpha,\beta) & \leq
\inf_{D}\lr{C_{\alpha,\beta}^{2} \rho_{D}^{2} + \varphi^{2}(s_{D}^{2})}
\end{align*}
Thus the bounds established in this study are sharp whenever
$\rho^{2}_{D} \asymp S_{D} v_{D}$, where $S_{D}^{2} = \sum_{j=1}^{D}
s_{j}^{-2}$, and $v_{D}^{2} = s_{D}^{-2}$, respectively. More
explicitly, if
$$
\sum_{j=1}^{D} \frac 1 {s_{j}^{4}} \asymp \frac 1 {s_{D}^{2}} \sum_{j=1}^{D} \frac 1 {s_{j}^{2}}.
$$
This concerns only the decay rate of the singular numbers $s_{j}$ of
the operator~$T$, and this holds for regularly varying singular
numbers, but this also holds true for $s_{j} \asymp \exp(- \gamma j),
j=1,2,\dots,$ thus covering severely ill-posed problems.
Remark that instead of the terms involved in (\ref{eq:rhod}), the quantities $S_{D}$ and $v_{D}$ have nice interpretation as strong
and weak variances of the spectral cut-off schemes.

\section{Relating the direct and inverse testing problem}
\label{sec:d+i}

For injective linear operators~$T$, the assertions "$f=0$" and "$T f
=0$" are equivalent. Hence, testing $H_{0}: f=0$ or testing $H_{0}: T
f =0$ is related to the same problem: we want to detect whether there
is signal in the data. Nevertheless, these testing problems are
different in the sense that the alternatives are not expressed in the
same way. Indeed, the inverse testing problem (considered in the
previous sections) corresponds to
\begin{equation}
H_0^I: f=0, \ \mathrm{against} \ H_1^I: f\in\mathcal{E}_\varphi, \ \| f\|^2 \geq (\rho^I)^2,
\label{eq:inverse_test}
\end{equation}
while the direct testing problem corresponds to test
\begin{equation}
H_0^D: T f=0,\ \mathrm{against} \ H_1^D: f\in\mathcal{E}_\varphi, \ \| T f\|^2 \geq (\rho^D)^2. 
\label{eq:direct_test}
\end{equation}
In this section, we investigate the similarities between these two
view points. In particular, we remark that both testing
problems are not equivalent in the sense that the alternatives do not
deal with the same object.

\subsection{Relating the separation rates}
\label{sec:sep-rates}

The authors in~\cite{MR2763215} discussed whether both problems are
related. The main result, Theorem~1, ibid. asserts that for a variety
of cases each minimax test~$\Phi_{\alpha}$ for the direct problem
($H_{0}: T f =0$) is also minimax for the related inverse problem
($H_{0}: f=0$). This fundamental results is based in Lemma~1,
ibid. Here we show that this lemma has its origin in
\emph{interpolation} in variable Hilbert scales, and we refer
to~\cite{MR2277542}. Actually we do not need the machinery as
developed there, but we may use the following special case, which may
directly be proved using Jensen's inequality.
\begin{lemma}[Interpolation inequality]\label{lem:interpol}
  Let $\varphi$ be from~(\ref{ellips}), and let $\Theta(u) := \sqrt u
  \varphi(u),\ u>0.$ If the function $u\mapsto
  \varphi^{2}\lr{\lr{\Theta^{2}}^{-1}(u)}$ is concave then
  \begin{equation}
    \norm{f}{} \leq \varphi\lr{\Theta^{-1}(\norm{T f}{})},\quad
    f\in\mathcal E_{\varphi}.
    \label{eq:interp}
  \end{equation}
\end{lemma}
The main result relating the direct and inverse testing problems is
the following.
\begin{thm}\label{thm:d+i}
  Let $\varphi$ be an index function with related function~$\Theta$,
  such that the function~$u\mapsto
  \varphi^{2}\lr{\lr{\Theta^{2}}^{-1}(u)}$ is concave.  Let
  $\Phi_{\alpha}$ be a level-$\alpha$ test for the direct problem
  $H_{0}^D: T f =0$ with uniform separation rate
  $\rho^D(\Phi_{\alpha},\mathcal E_{\Theta},\beta)$.  Then
  $\Phi_{\alpha}$ constitutes a level-$\alpha$ test for the inverse
  problem $H_{0}^I: f=0$ with uniform separation rate
  $$
  \rho^I(\Phi_{\alpha},\mathcal E_{\varphi},\beta) \leq
  \varphi\lr{\Theta^{-1}\lr{\rho^D(\Phi_{\alpha},\mathcal
      E_{\Theta},\beta)}}.
$$
Consequently we have for the minimax separation rates that
\begin{equation}
  \label{eq:minimax-rates-compare}
  \rho^I(\mathcal E_{\varphi},\alpha,\beta) \leq
  \varphi\lr{\Theta^{-1}\lr{\rho^D(\mathcal E_{\Theta},\alpha,\beta)}}.
\end{equation}
\end{thm}
\begin{proof}
  Clearly, the test $\Phi_{\alpha}$ is a level-$\alpha$ test for both
  problems, and we need to control the second kind error. But if
  $\norm{f}{} \geq
  \varphi\lr{\Theta^{-1}\lr{\rho^D(\Phi_{\alpha},\mathcal
      E_{\Theta},\beta)}}$ then Lemma~\ref{lem:interpol} yields that
  $\norm{T f}{} \geq \rho^D(\Phi_{\alpha},\mathcal E_{\Theta},\beta)$,
  and the assertion is a consequence of the properties of the test for
  the direct problem.

  If $\Phi_{\alpha}$ was minimax for the direct problem then the
  corresponding minimax rate for the inverse problem must be dominated
  by $\varphi\lr{\Theta^{-1}\lr{\rho^D(\mathcal
      E_{\Theta},\alpha,\beta)}}$, which
  gives~(\ref{eq:minimax-rates-compare}).
\end{proof}
\begin{remark}
  In many cases the bound~(\ref{eq:minimax-rates-compare}) actually is
  an asymptotic equivalence
  \begin{equation}
    \label{eq:minimac-compare-equivalence}
    \varphi^{-1}\lr{\rho^I(\mathcal E_{\varphi},\alpha,\beta) } \asymp 
    \Theta^{-1}\lr{\rho^D(\mathcal E_{\Theta},\alpha,\beta) },\quad
    \sigma\to 0.
  \end{equation}
  It may be enlightening to see this on the base of
  Example~\ref{xmpl:moderate}. Recall that the function~$\varphi$ was
  given as $\varphi(u) = u^{s/(2t)}$. The corresponding rate is known
  to be minimax, and we obtain that
$$
\varphi^{-1}\lr{\rho^I(\mathcal E_{\varphi},\alpha,\beta) } \asymp
\sigma^{\frac{4t}{2s + 2t + 1/2}}.
$$
We turn to the direct problem, for which the corresponding smoothness
class is $\mathcal E_{\Theta}$ for the function~$\Theta(u) = u^{2/(2t)
  + 1/2} = u^{(s+t)/(st)}$. This corresponds to $\mu = s+t$
in~\cite[Tbl.~2]{MR2763215}, yielding the separation rate
$\rho(\mathcal E_{\Theta},\alpha,\beta) \asymp \sigma^{2(s+t)/(2s + 2t
  + 1/2)}$, which in turn gives
$$
\Theta^{-1}\lr{\rho^D(\mathcal E_{\Theta},\alpha,\beta) } \asymp
\sigma^{\frac{4t}{2s + 2t + 1/2}},
$$
and hence~(\ref{eq:minimac-compare-equivalence}) for moderately
ill-posed problems.

Similarly, this holds for severely ill-posed problems, and we omit
details.
\end{remark}
We emphasize that, by virtue of Theorem~\ref{thm:d+i}, any lower bound
for the minimax separation rate in the inverse testing problem yields
a lower bound for the corresponding direct problem.

\begin{remark}
Thanks to Theorem \ref{thm:d+i}, it is possible to prove that in all
the cases considered in this paper, a test minimax for
(\ref{eq:direct_test}) will be also minimax for
(\ref{eq:inverse_test}). Nevertheless, the reverse is not true. We
will not dwell into details, instead we refer to \cite{MR2763215} for
a detailed discussion on this subject.  
\end{remark}

\subsection{Designing tests for the direct problem}
\label{sec:design-tests}

The coincidence in~(\ref{eq:minimac-compare-equivalence}) is not by
chance and we indicate a further result in this direction. Recall
from~\ref{thm:main-linear} that the value of $\tau_{\ast}=
\tau_{\ast}^{\IP}$ was obtained from~(\ref{eq:tauast}), and hence that
we actually have $\rho(\mathcal E_{\varphi},\alpha,\beta) \asymp
\varphi(\tau_{\ast}^{\IP})$, such that the left hand side
in~(\ref{eq:minimac-compare-equivalence}) equals
$\tau_{\ast}^{\IP}$. We shall see next that the corresponding value
$\tau_{\ast}= \tau_{\ast}^{\DP}$ is obtained from the same
equation~(\ref{eq:tauast}) when basing the direct test on the family
$\widehat{T R_{\tau}} = T R_{\tau}$ with family $R_{\tau} =
\ga(T^{\ast} T)T^{\ast}$ as in \S~\ref{sec:linreg}.  Then $T R_{\tau}
= \ga(T T^{\ast}) T T^{\ast}$, and we bound its variance and weak
variance, next.
\begin{lemma}\label{lem:sv-direct}
  Let $\tilde R_\tau=T R_{\tau} = \ga(T T^{\ast}) TT^{\ast}$ and
  denote by resp. $\tilde S_\tau^2$ and $\tilde v_\tau^2$ the
  corresponding strong and weak variance. If Assumption~\ref{ass:T-HS}
  holds then
  \begin{enumerate}
  \item $\tilde S_{\tau}^{2} \leq (\gamma_{0} +
    \gamma_{\ast})\gamma_{0} \sigma^{2} \cN(\tau),\ \tau>0$, and
  \item $\tilde v_{\tau}^{2} \leq \sigma^{2} \gamma_{0}^{2}$.
  \end{enumerate}
\end{lemma}
We also need to bound the bias $\norm{T f - T f_{\tau}}{}$ with
$f_{\tau} = \ga(T^{\ast} T)T^{\ast} T f$
\begin{lemma}
  Assume that $f\in\mathcal{E}_\varphi$. If the regularization $\ga$
  has qualification $\Theta$ with constant $\gamma $ then
$$
\norm{T f - T f_{\tau}}{} \leq \gamma \Theta(\tau).
$$ 
\end{lemma}
\begin{proof}
  Since $f_\tau=R_\tau T f$, we get that
$$ \| T f-T f_\tau \| = \| T f - g_\tau(T^*T)TT*T f \| = \| r_\tau(T^*T) T f \| ,$$
which is bounded by $\gamma \Theta(\tau)$ as soon as $f\in
\mathcal{E}_\varphi$ and $g_\tau$ has qualification $\Theta$.
\end{proof}

We recall from~\S~\ref{sec:sepradius} the
quantity~$r^{2}(\Phi_{\alpha},\beta):= C_{\alpha,\beta}
S_{\tau}v_{\tau}$, where we now consider $\tilde R_{\tau}$ and $\tilde
v_{\tau}$ from Lemma~\ref{lem:sv-direct} for bounding $\norm{T
  f}{}^{2} \geq C_{\alpha,\beta} \tilde S_{\tau} \tilde v_{\tau}$ from
below.
\begin{corollary}
  Suppose that $\ga$ is a regularization which has
  qualification~$\Theta$, $f\in\mathcal E_{\varphi}$ and that
  Assumption~\ref{ass:T-HS} holds.  Let $\tau_{\ast}^{DP}$ be chosen
  from the equation
  \begin{equation}
    \label{eq:tauastDP}
    \sigma ^{2} = \frac{\Theta^{2}(\tau)}{\sqrt{\cN(\tau)}}.    
  \end{equation}
  Then
$$
\inf_{\tau>0} \lr{r^{2}(\Phi_{\alpha},\beta) + \norm{T f - T
    f_{\tau}}{}^{2}} \leq \lr{C_{\alpha,\beta} \sqrt{(\gamma_{0}
    +\gamma_{\ast})\gamma_{0}}\gamma_{0} +
  \gamma^{2}}\Theta^{2}(\tau_{\ast}^{\DP}).
$$
\end{corollary}
We stress that the equation~(\ref{eq:tauastDP}) for determining
$\tau_{\ast}^{\DP}$ is the same equation as~(\ref{eq:tauast}), since
$\Theta^{2}(\tau) = \tau \varphi^{2}(\tau)$, and this explains the
identical asymptotics in~(\ref{eq:minimac-compare-equivalence}) as
being equal to $\tau_{\ast}^{\DP} = \tau_{\ast}^{\IP}$.

This result sheds light to another interesting problem: If we want to
use the regularization $T R_{\tau}$, and if we want to have this
optimal performance properties then the underlying
regularization~$\ga$ must have \emph{higher qualification}~$\Theta$
for the direct problem as compared for its use in inverse testing
requiring qualification~$\varphi$, only. This cannot be seen when
confining to spectral cut-off, but this problem is relevant when
considering other regularization schemes for testing. It is thus
interesting to design estimators for $g= T f$ which do not rely on
estimation of $f$. However, since the data $Y$ do not belong to the
space~$K$ either discretization or some other kind of preconditioning
is necessary in order to estimate $g=T f$ from the data~$Y$.  Such
direct estimation is simple by using projection schemes, and we
exhibit the calculus for one-sided discretization. As
in~\S~\ref{sec:projschemes} we choose finite ($m$) dimensional
subspaces~$Y_{m}\subset K$, with corresponding projections~$Q_{m}$ and
consider the data
$$
Q_{m}Y = Q_{m}g + \sigma Q_{m}\xi,\quad m\in\nat.
$$
This approach is called \emph{dual least squares} scheme in
regularization, see~\cite{MR0488721}.  Here it is easy to see that
$S_{m}^{2} = \tr{Q_{m}^{\ast}Q_{m}} =m$, while $v_{m}^{2} =
\norm{Q_{m}}{}^{2}=1$.  In order to continue we just need that the
chosen projections have degree of approximation~$\Theta$, i.e,\ there
is $C_{D}$ for which $\norm{(I - Q_{m})\Theta(T T^{\ast})}{}\leq C_{D}
\Theta(s_{m+1}^{2}),\ m=1,2\dots.$ With this requirement at hand we
can continue as if the projections~$Q_{m}$ were the projections onto
the first $m$ singular elements in the svd of $T$.  In particular we
have the upper bound on the separation radius
$$
\rho(\mathcal E_{\Theta},\alpha,\beta) \leq \max\set{\Dab,C_{D}^{2}}
\inf_{m} \lr{ \sigma^{2} \sqrt m + \Theta^{2}(s_{m+1}^{2})},
$$
similar to corresponding results obtained for spectral cut-off
in~\cite{MR1935648,MR2763215}, and we omit further details.

\section{Adaptation to the smoothness of the alternative}
\label{s:adaptation}
It seems clear from Section 4 that the optimality of the considered
tests strongly depends on the regularity (smoothness) of the
alternative. In this section, we propose data-driven tests that
automatically adapt to the unknown smoothness parameter.  The
adaptation issue in test theory has widely been investigated. For more
details on the subject, we refer for instance to \cite{BHL_2003},
\cite{Spokoiny_1996} in the direct setting (i.e. $T=Id$) or
\cite{ISS_2012} in the inverse case for an adaptive scheme based on
the singular value decomposition of the operator.

First, we propose a general adaptive scheme. Then, we apply this
approach to linear regularization over ellipsoids. This methodology
can also be extended to projection schemes. For the sake of brevity,
this extension is not discussed here.

\subsection{A general scheme for adaptation}
Assume that we have at our disposal a finite collection
$(R)_{R\in\mathcal{R}}$ of regularization operators satisfying
Assumption~\ref{ass:T-HS}. Then, we can associate to each operator $R$
a level-$\alpha$ test $\Phi_{\alpha,R}$. Our aim in this section is to
construct a test that mimics the behavior of the best possible test
among the family $\mathcal{R}$.  Let $|\mathcal{R}|$ denotes the
cardinality of the family $\mathcal{R}$. We define our adaptive test
$\Phi_\alpha^\star$ as
\begin{equation}
  \Phi_\alpha^\star = \max_{R\in\mathcal{R}} \Phi_{\frac{\alpha}{|\mathcal{R}|},R}.
  \label{eq:adaptive_test}
\end{equation}
The performance of $\Phi_\alpha^\star$ is summarized in the following
proposition.

\begin{proposition}
  \label{prop:adaptation}
  The test introduced in (\ref{eq:adaptive_test}) is a level-$\alpha$
  test. Moreover
$$ P_f(\Phi_\alpha^\star =0) \leq \beta,$$ 
as soon as
$$ \| f\|^2 \geq 2 \inf_{R\in \mathcal{R}}
\lr{r^{2}(\Phi_{\frac{\alpha}{|\mathcal{R}|},R},\beta) + \norm{f -
    f_{R}}{}^{2}},$$ where the term $r^{2}$ has been introduced in
(\ref{eq:r-def}).
\end{proposition}

\begin{proof}
  We first remark that
  \begin{align*}
    P_{H_0}(\Phi_\alpha^\star =1)
    & =  P_{H_0}\left( \max_{R\in\mathcal{R}} \Phi_{\frac{\alpha}{|\mathcal{R}|},R} =1 \right),\\
    & =  P_{H_0}\left( \bigcup_{R\in\mathcal{R}} \Phi_{\frac{\alpha}{|\mathcal{R}|},R} =1 \right),\\
    & \leq \sum_{R\in\mathcal{R}} P_{H_0}\left(
      \Phi_{\frac{\alpha}{|\mathcal{R}|},R} =1 \right) = \alpha,
  \end{align*}
  since $P_{H_0}(\Phi_{\frac{\alpha}{|\mathcal{R}|},R} =1 ) =
  \alpha/|\mathcal{R}|$ for all $R\in \mathcal{R}$. Hence,
  $\Phi_\alpha^\star$ is a level-$\alpha$ test. Now, we can
  investigate the second kind error. Using simple algebra, we get that
  \begin{align*}
    P_f(\Phi_\alpha^\star =0)
    & =  P_{H_0}\left( \max_{R\in\mathcal{R}} \Phi_{\frac{\alpha}{|\mathcal{R}|},R} =0 \right),\\
    & =  P_{H_0}\left( \bigcap_{R\in\mathcal{R}} \Phi_{\frac{\alpha}{|\mathcal{R}|},R} =0 \right),\\
    & \leq \inf_{R\in\mathcal{R}} P_{H_0}\left(
      \Phi_{\frac{\alpha}{|\mathcal{R}|},R} =0 \right).
  \end{align*}
  We can conclude using (\ref{eq:cond_f}). 
\end{proof}

Proposition \ref{prop:adaptation} proves that the detection radius
associated to the test defined in (\ref{eq:adaptive_test}) is close to
the smallest possible one among the family $\mathcal{R}$. Thus, we
must design the set $\mathcal R$ according to two requirements. First,
the cardinality $\abs{\mathcal R}$ should be small, in order not to
enlarge the detection radius too much.  Indeed, the following holds
true.
\begin{lemma}\label{lem:cardinalityR}
  Let $C_{\alpha,\beta}^*$ the term introduced
  in~(\ref{eq:dabast}). If the family $\mathcal R$ of regularization
  schemes has cardinality $M:= \abs{\mathcal R}\geq 1$, then
$$
C^*_{\alpha/M,\beta} \leq C_{\alpha,\beta}^* + 2\sqrt{2 \log(M)}.
$$
If $M \geq 4$ then $C^*_{\alpha/M,\beta} \leq
(C^*_{\alpha,\beta}+2\sqrt{2})\sqrt{\log(M)}$.
\end{lemma}
\begin{proof}
  We first observe that $x_{\alpha/M} = x_{\alpha} +
  \log(M)$. Therefore we conclude that
  \begin{align*}
    C_{\alpha,\beta}^{\ast}& = 2 \sqrt{x_{\beta}} + 2 \sqrt{2
      x_{\alpha/M}} ,\\
    & =  C_{\alpha,\beta}^{\ast} + 2 \lr{\sqrt{2 x_{\alpha} + 2\log(M)} -\sqrt{2 x_{\alpha}} }.\\
    & \leq C_{\alpha,\beta}^{\ast} + 2 \sqrt{2 \log(M)} .
  \end{align*}
  The second assertion is trivial, because $\log(M)>1$ for $M\geq 4$.
\end{proof}
Therefore, the price to pay for using $\Phi_\alpha^\star$ is a term of
order $\sqrt{\log(\abs{\mathcal{R}})}$, up to some condition on the
behavior of the effective dimension (see Theorem \ref{thm:summarize}
below).  On the other hand, the set $\mathcal R$ should be rich enough
to keep the detection radius on the size of the best possible bound,
as such was established Theorems~\ref{thm:main-linear} and
\ref{thm:main-projection}.

In the following, we propose practical situations where such an
adaptive scheme can be used. In particular, we propose families of
regularizations operators with controlled size and prove that the
adaptive test~$\Phi_{\alpha}^{\ast}$ attains the minimax rate of
testing (up to a $\log\log$ term) for a proper choice of $\mathcal R$.

\begin{remark}
  In the test (\ref{eq:adaptive_test}), each regularization operator
  $R\in\mathcal{R}$ is associated to a test
  $\Phi_{\frac{\alpha}{|\mathcal{R}|},R}$ having the same level
  $\alpha/|\mathcal{R}|$. It is nevertheless possible to use more
  refined approaches, leading to an improvement of the power of the
  test (in terms of the constants). We refer to
  \cite[Eq.~(2.2)]{FL_2006}, however in a slightly  different setting. 
\end{remark}

\subsection{Application to linear regularization}

We will exhibit the use of the general methodology for tests based on
linear regularization.

Let $g_\tau$ be  a given regularization. We associate to each function
$g_\tau$ the operator $R_\tau$ and we deal with the family
$\mathcal{R}=(R_\tau)_{\tau>0}$. In order to apply
Proposition~\ref{prop:adaptation} we need to specify a finite subset
$\mathcal R \subset (0,\infty)$ on which the test
$\Phi_{\alpha}^{\ast}$ will be based on.  To this end we will use an
exponential grid.  Given an initial value $\tau_{\max}$, and a tuning
parameter $0 < q < 1$ we consider the exponential grid
\begin{equation}
  \label{eq:grid}
  \Delta_{q} := \set{\tau = q^{j}\tau_{\max},\quad j=0,\dots, M-1},\quad
  \text{for some } M>1.
\end{equation}
Then we use the adaptive test
\begin{equation}
  \Phi_{\alpha}^{\star} = \max_{\tau\in\Delta_{q}} \Phi_{\alpha/M,\tau}.
  \label{eq:adaptive_test2}
\end{equation}
The result from Proposition~\ref{prop:adaptation} can be rephrased as
follows.  By virtue of Lemma~\ref{lem:cardinalityR}, and using the
bounds from Lemma~\ref{lem:Svbounds} \& Proposition~\ref{pro:quali},
respectively, we find that the test $\Phi_{\alpha}^{\ast}$ bounds the
error of the second kind by $\beta$ as soon as
$$
\norm{f^{2}}{} \geq C(\alpha,\beta) \inf_{\tau\in\Delta_{q}}
\lr{\sqrt{\log(M)}\sigma^{2} \frac{\sqrt{\cN(\tau)}}{\tau}+ \log(M)
  \frac{\sigma^{2}} {\tau} + \varphi^{2}(\tau)},
$$
for some explicit constant $C(\alpha,\beta)$. We shall now show, how
we can specify the numbers~$0 < \tau_{\min} < \tau_{\max}$ such that
this is of the order of the separation radius (up to a
$\log\log$-factor).

The cardinality~$M$ obeys~$\tau_{\min}:= q^{M-1}\tau_{\max}$, and
hence $M := \log_{\frac 1 q }(\tau_{\max}/\tau_{\min})$.  Obviously we
have that
\begin{eqnarray*}
  \lefteqn{\inf_{\tau_{\min} \leq \tau \leq \tau_{\max}} \lr{\sqrt{\log(M)}\sigma^{2} \frac{\sqrt{\cN(\tau)}}{\tau} 
      + \log(M)   \frac{\sigma^{2}}{\tau} +   \varphi^{2}(\tau)}}\\  
  & \leq & \inf_{\tau\in\Delta_{q} }\lr{\sqrt{\log(M)}\sigma^{2} \frac{\sqrt{\cN(\tau)}}{\tau} + \log(M)   \frac{\sigma^{2}}{\tau} +
    \varphi^{2}(\tau)}
\end{eqnarray*}
The reverse is also true (up to some constant), as proved in the
following lemma.
\begin{lemma}[{cf.~\cite[Proof of Thm.~3.1]{JM2012}}]
  \label{lem:all2grid}
  We have that
  \begin{eqnarray*}
    \label{eq:all2grid}
    \lefteqn{\inf_{\tau_{\min}\leq  \tau\leq \tau_{\max}} \lr{\sqrt{\log(M)}\sigma^{2} \frac{\sqrt{\cN(\tau)}}{\tau}+ \log(M)   \frac{\sigma^{2}}{\tau} + \varphi^{2}(\tau)}}\\
    & \geq & q^{3/2}  \inf_{\tau\in\Delta_{q} }\lr{\sqrt{\log(M)}\sigma^{2} \frac{\sqrt{\cN(\tau)}}{\tau} 
      + \log(M)   \frac{\sigma^{2}}{\tau} +\varphi^{2}(\tau)}. \nonumber
  \end{eqnarray*}
\end{lemma}
\begin{proof}
  For any $\tau$ with $\tau_{\min} < \tau \leq \tau_{\max}$ we find an
  index $1 \leq j \leq M$ for which $\tau_{j} < \tau \leq
  \tau_{j}/q$. The crucial observation is that the function~$\tau \to
  \frac{\sqrt{\cN(\tau)}}{\tau}$ is decreasing, whereas the function
  $\tau \to \sqrt{\tau\cN(\tau)} = \tau^{3/2}
  \frac{\sqrt{\cN(\tau)}}{\tau}$ is increasing, which can be seen from
  spectral calculus.  Therefore, by using the above monotonicity we
  see that
  \begin{align*}
    \lefteqn{\sqrt{\log(M)}\sigma^{2} \frac{\sqrt{\cN(\tau)}}{\tau} + \log(M)   \frac{\sigma^{2}}{\tau} + \varphi^{2}(\tau) }\\
    & \geq
    \sqrt{\log(M)}\sigma^{2} \frac{\sqrt{\cN(\tau_{j}/q)}}{\tau_{j}/q} + \log(M)   \frac{\sigma^{2}}{\tau_j/q} + \varphi^{2}(\tau_{j}) \\
    & = \sqrt{\log(M)}\sigma^{2} \lr{\frac{\tau_{j}}{q}}^{-3/2}
    \lr{\frac{\tau_{j}}{q}}^{3/2}
    \frac{\sqrt{\cN(\tau_{j}/q)}}{\tau_{j}/q} + q \log(M)   \frac{\sigma^{2}}{\tau_j}+ \varphi^{2}(\tau_{j}) \\
    &\geq \sqrt{\log(M)}\sigma^{2} \lr{\frac{\tau_{j}}{q}}^{-3/2}
    {\tau_{j}}^{3/2}
    \frac{\sqrt{\cN(\tau_{j})}}{\tau_{j}} + q^{3/2} \log(M)   \frac{\sigma^{2}}{\tau_j} + \varphi^{2}(\tau_{j}) \\
    & \geq q^{3/2}
    \lr{\sqrt{\log(M)}\sigma^{2}\frac{\sqrt{\cN(\tau_{j})}}{\tau_{j}}
      + \log(M) \frac{\sigma^{2}}{\tau_j}+ \varphi^{2}(\tau_{j})},
  \end{align*}
  from which the proof can easily be completed.
\end{proof}
We shall next discuss the choices of $\tau_{min}$ and $\tau_{\max}$.
First, the natural domain of definition of the smoothness
function~$\varphi$ is on $(0,\norm{T^{\ast}T}{}]$, such that the
choice $\tau_{\max} = \norm{T^{\ast}T}{}$ is natural. In this case the
size of $\sqrt{\log(M)}\sigma^{2}\frac{\sqrt{\cN(\tau)}}{\tau} +
\log(M) \frac{\sigma^{2}}{\tau} + \varphi^{2}(\tau_{\max})$ is at
least $\varphi^{2}(\norm{T^{\ast}T}{})$ no matter how small the noise
level~$\sigma$ was.  The next result indicates that we can find
$\tau_{\min}$ in such a way that we can remove the restriction to
$\tau > \tau_{\min}$ if there is some 'minimal' smoothness in the
alternative.
\begin{lemma}
  \label{lem:minsmooth}
  Let $\tau_{\min} = \tau_{\min}(M)$ satisfy
  \begin{equation}
    \label{eq:taumin}
    \sqrt{\log(M)}\sigma^{2}
    \frac{\sqrt{\cN(\tau_{\min})}}{\tau_{\min}} \geq 1.
  \end{equation} If the smoothness $\varphi$ obeys
  $\varphi(\tau_{\min})\leq 1$ then for $0< \tau \leq \tau_{\min}$ we
  have that
  \begin{multline*}
    \lefteqn{\log(M)\sigma^{2} \frac{\sqrt{\cN(\tau)}}{\tau} + \log(M)
      \frac{\sigma^{2}}{\tau} + \varphi^{2}(\tau) }\\
    \geq \frac 1 2 \lr{ \sqrt{\log(M)}\sigma^{2}
      \frac{\sqrt{\cN(\tau_{\min})}}{\tau_{\min}} + \log(M)
      \frac{\sigma^{2}}{\tau_{min}} + \varphi^{2}(\tau_{\min}) }.
  \end{multline*}
\end{lemma}
\begin{proof}
  For $\tau < \tau_{\min}$ this easily follows from
  \begin{align*}
    \sqrt{\log(M)}\sigma^{2} \frac{\sqrt{\cN(\tau)}}{\tau} +
    \varphi^{2}(\tau) &\geq
    \sqrt{\log(M)}\sigma^{2} \frac{\sqrt{\cN(\tau_{\min})}}{\tau_{\min}} \geq 1 \\
    &\geq \frac 1 2 \lr{\sqrt{\log(M)}\sigma^{2}
      \frac{\sqrt{\cN(\tau_{\min})}}{\tau_{\min}} +
      \varphi^{2}(\tau_{\min})},
  \end{align*}
  which proves the assertion.
\end{proof}
\begin{remark}
  For given $\sigma>0$ the condition from~(\ref{eq:taumin}) can always
  be satisfied.  Below we shall further specify this as follows. If
  $\tau_{\max}$ is chosen as $\norm{T^{\ast}T}{}$ then $\cN(\tau_{\max})\geq
  1/2$, such that
$$
\sqrt{ \log(M)}\sigma^{2} \frac{\sqrt{\cN(\tau_{\min})}}{\tau_{\min}}
\geq \frac{\sqrt{\log(M)}\sigma^{2}}{\sqrt 2\tau_{\min}} =
\frac{\sqrt{\log(M)}\sigma^{2} q^{1-M}}{\sqrt 2 \norm{T^{\ast}T}{}}
\geq\frac{1}{\sqrt 2 \norm{T^{\ast}T}{}} \frac{\sigma^{2}}{q^{M}}.
$$
Thus the condition~(\ref{eq:taumin}) holds for
$$
M \geq \log_{1/q} \lr{\sqrt 2\norm{T^{\ast}T}{} } + \log_{1/q}
(1/\sigma^{2}).
$$
\end{remark}
We summarize the above considerations.
\begin{proposition}
  \label{prop:summarize}
  Suppose that $M$ and $\tau_{\min}$ are chosen such
  that~(\ref{eq:taumin}) holds. If the smoothness function~$\varphi$
  obeys $\varphi(\tau_{\min}) \leq 1$ then
  \begin{align*}
    \label{eq:summarize}
    \lefteqn{\inf_{\tau \in\Delta_{q}} \lr{\sqrt{\log(M)}\sigma^{2}
        \frac{\sqrt{\cN(\tau)}}{\tau} + \log(M)
        \frac{\sigma^{2}}{\tau} + \varphi^{2}(\tau)}}\\
    & \leq q^{-3/2} 2 \inf_{0< \tau \leq \tau_{\max}}
    \lr{\sqrt{\log(M)}\sigma^{2} \frac{\sqrt{\cN(\tau)}}{\tau} +
      \log(M) \frac{\sigma^{2}}{\tau} + \varphi^{2}(\tau)}.
  \end{align*}
\end{proposition}
The following result summarizes the above considerations; it asserts
that the test $\Phi_{\alpha}^{\star}$ appears to be minimax (up to a
$\log\log$ term) in many cases.
\begin{thm}\label{thm:summarize}
  Let $\alpha$, $\beta$ be fixed and $\Phi_{\alpha}^{\star}$ the test
  defined in (\ref{eq:adaptive_test2}).  Suppose that $\tau_{\max}=
  \norm{T^{\ast}T}{}$, $\tau_{\min}$ is chosen such that $M \geq
  \log_{1/q} \lr{\sqrt 2\norm{T^{\ast}T}{} } + \log_{1/q}
  (1/\sigma^{2})$. Let $\tau_{\ast}$ be given from
  \begin{equation}
    \label{eq:tauastadapt}
    \varphi^{2}(\tau_{\ast}) = \sigma ^{2} \sqrt{\log\log_{1/q}(\frac 1
      {\sigma^{2}})} \frac{\sqrt{\cN(\tau_{\ast})}}{\tau_{\ast}},   
  \end{equation}
  If the underlying smoothness obeys $\varphi(\tau_{\min})\leq 1$ and
  if
  \begin{equation}
    \label{eq:cond5.5}
    \frac{\log\log(\frac 1 {\sigma^{2}})}{\cN(\tau_{\ast})} =\littleo(1) \quad \text{as} \quad\sigma\rightarrow 0,
  \end{equation}
  then there is a constant $C>0$ such that
$$ \rho^{2}(\Phi_\alpha^\star,\beta,\mathcal{E}_\varphi) \leq C 
\inf_{0< \tau \leq \tau_{\max}}
\lr{\sigma^{2}\sqrt{\log\log_{1/q}(\frac 1 {\sigma^{2}})}
  \frac{\sqrt{\cN(\tau)}}{\tau} + \varphi^{2}(\tau)}.
$$
In particular, as $\sigma\searrow 0$ we have that $\tau_{\ast}\searrow
0$, and hence that there is a constant~$D = D(\alpha,\beta)$ such that
$$
\rho(\Phi_\alpha^\star,\beta,\mathcal{E}_\varphi) \leq D
\varphi(\tau_{\ast}),\quad \text{as } \sigma\searrow 0.
$$
\end{thm}
We shall indicate that the assumption~(\ref{eq:cond5.5}) is valid in
many cases
.
\begin{lemma}
  If there is a constant $c>0$ such that the effective dimension obeys
  \begin{equation}
    \label{eq:cN-lb}
    \cN(\tau) \geq c \log(1/\tau),
  \end{equation}
  and if the smoothness increases at least as
  \begin{equation}
    \label{eq:varphi-ub}
    \varphi(\tau) \leq \lr{\log\log_{1/q}(1/\tau)}^{4},
  \end{equation}
  as $\tau \to 0$, then~(\ref{eq:cond5.5}) is valid.
\end{lemma}
\begin{proof}
  The parameter~$\tau_{\ast}$ is determined
  from~(\ref{eq:tauastadapt}), and under~(\ref{eq:varphi-ub}) we find
  that
$$
\sigma ^{4} {\log\log_{1/q}(\frac 1 {\sigma^{2}})} =
\frac{\tau_{\ast}^{2} \varphi^{4}(\tau_{\ast})}{\cN(\tau_{\ast})} \leq
\frac{\tau_{\ast}^{2}\log\log_{1/q}(1/\tau_{\ast})}{\cN(\tau_{\ast})}
\leq\tau_{\ast}^{2}\log\log_{1/q}(\frac 1{\tau_{\ast}}),
$$
provided that $\tau_{\ast}$ is small enough.  Monotonicity yields that
$\sigma^{2}\leq \tau_{\ast}$. But then $\log\log(1/\sigma^{2}) \leq
\log\log(1/\tau_{\ast})$, and we conclude that
$$
\frac{\log\log(1 /{\sigma^{2}})}{\cN(\tau_{\ast})}\leq
\frac{\log\log(1/\tau_{\ast})}{\cN(\tau_{\ast})} \leq \frac 1 c
\frac{\log\log(1/\tau_{\ast})}{\log(1/\tau_{\ast})} = \littleo(1),
$$
as $\sigma$, and hence $\tau_{\ast}$, tend to zero.
\end{proof}
\begin{remark}
  This result covers many of the interesting cases, in particular the
  ones from Examples~\ref{xmpl:moderate} \& \ref{xmpl:severe}.  In
  these cases Theorem~\ref{thm:summarize} exhibits that the separation
  radii obey
  \begin{align*}
    \rho(\Phi^{\ast}_{\alpha,\tau_*},\beta,\mathcal{E}_\varphi) & \leq
    D \lr{\sigma^2\sqrt{\log\log \frac 1 {\sigma^{2}}}}^{s/(2s + 2t +
      1/2)}, \
    \text{ and }\\
    \rho(\Phi^{\ast}_{\alpha,\tau_*},\beta,\mathcal{E}_\varphi) & \leq
    D\log^{-s}(1/\sigma^{2}),
  \end{align*}
  respectively. In particular we see that adaptation does not pay an
  additional price for severely ill-posed testing problems.
\end{remark}

\begin{remark}
  A similar approach can be used when basing the adaptive test on a
  family of projection schemes. In this case we use a finite family of
  dimensions
$$
\Delta_{2,j_{0}} := \set{m = 2^{j + j_{0}},\quad j=0,\dots, M-1},
$$
and consider projection schemes with spaces $X_{m},\ Y_{n(m)}$ for $m
\in \Delta_{2,j_{0}}$. The above reasoning applies, taking into
account the correspondence between regularization parameter $\tau$ in
linear regularization schemes, and dimensions~$m \sim 1/\tau$ . For
the sake of brevity, this will not be discussed in this paper.
\end{remark}

 \appendix 
\section{Inequalities for Gaussian elements in Hilbert space}

\begin{lemma}
  \label{TL}
  Let $X$ a Gaussian random variable having values in $H$. Then, for
  all $x>0$,
$$ 
P\lr{\norm{ X}{}^{2} - \expect\norm{X}{}^{2} \geq x^2 + 2x
  \sqrt{\expect\norm{X}{}^{2}}} \leq \exp\lr{  -\frac{x^2}{2v^2}}.
$$
where
$$
v^{2}:= \sup_{\|\omega\|\leq 1} \expect |\langle X,\omega
\rangle|^2.
$$
\end{lemma}
\begin{proof}
  Using, the Cauchy-Schwarz inequality, we first observe  that
$$ 
\left( \expect[\| X \|] +x \right)^2 
\leq \expect\| X \|^2 + x^2 + 2x \sqrt{\expect\| X \|^2}.
$$
Hence, we get
\begin{align*}
  P\lr{\norm{ X}{}^{2} - \expect\norm{X}{}^{2} \geq x^2 + 2x
  \sqrt{\expect\norm{X}{}^{2}}}
  & \leq P\left(\| X \|^2 \geq  \left( \expect[\| X \|] +x \right)^2 \right),\\
  & =  P\left(\| X \| \geq  \expect[\| X \|] +x  \right),\\
  & \leq  \exp \left( -\frac{x^2}{2v^2} \right),
\end{align*}
where for the last inequality we have used~\cite[Lemma~3.1]{MR1102015}.
\end{proof}

\begin{lemma}
  \label{lem:deviation_left}
  Let $R Y$ be as in~(\ref{eq:estimator}). Then
$$ P_f \left( \| R Y \|^2 - \expect_f \| R Y \|^2  \leq - 2
  \sqrt{\Sigma x_{\beta}} \right) \leq \beta,$$ where $\Sigma$ is
from~(\ref{eq:sigma}).
\end{lemma}
\begin{proof}
The proof is a direct extension of the one proposed in \cite{LLM_2012}
for a spectral cut-off approach. 
\end{proof}

\bibliography{inverse_test} 
\bibliographystyle{plain}

\end{document}